\documentclass[11pt, reqno]{amsart}
 
 \usepackage{amsmath}
 \usepackage{amssymb}
\usepackage{bm}

\DeclareMathAccent{\mathring}{\mathalpha}{operators}{"17}

 \usepackage{color}

\newcommand{\mysection}[1]{\section{#1}
      \setcounter{equation}{0}}

\newtheorem{theorem}{Theorem}[section]
\newtheorem{lemma}[theorem]{Lemma}

\newtheorem{corollary}[theorem]{Corollary} 

\theoremstyle{definition}
\newtheorem{assumption}{Assumption}[section]

\theoremstyle{remark}
\newtheorem{remark}{Remark}[section]

\newcommand{\tr}{\text{\rm tr}\,}
\newcommand{\loc}{\text{\rm loc}}

 \makeatletter 
 \def\dashint{%
 \operatorname%
 {\,\,\text{\bf--}\kern-.98em\DOTSI\intop\ilimits@\!\!}}
\def\ninf{\qopname\relax\@empty{inf\phantom{p}\!\!\!}}
 \makeatother

\newcommand\gb{\mathfrak{b}}

\newcommand\bbeta{\text{\raise-.2ex\hbox{$\bm{\beta}$}}}

\newcommand\esssup{\operatornamewithlimits{esssup}}

\newcommand\bR{\mathbb{R}}

\newcommand\bS{\mathbb{S}}

\newcommand\cB{\mathcal{B}}

\newcommand\cF{\mathcal{F}}

\newcommand\frM{\mathfrak{M}}

\newcommand\dist{{\rm dist}\,}

\newcommand\Vol{\text{\rm Vol}\,}

\begin{document}

\title[Stochastic processes with drift in $L_{d}$]
{On stochastic It\^o processes with drift in $L_{d}$}

\author{N.V. Krylov}
 
\email{nkrylov@umn.edu}
\address{127 Vincent Hall, University of Minnesota,
 Minneapolis, MN, 55455}



\begin{abstract}
   For It\^o stochastic processes
in $\bR^{d}$ with drift in $L_{d}$ 
Aleksandrov's type estimates are established
in the elliptic and parabolic settings.
They are applied to estimating the resolvent
operators of the corresponding elliptic
and parabolic operators in $L_{p}$
and $L_{p+1}$, respectively, where $p\geq d$.

\end{abstract}

\maketitle

\mysection{Introduction}

Let  $d_{1}$ be an integer $>1$,
$(\Omega,\cF,P)$ be a complete probability space,
and let $(w_{t},\cF_{t})$ be a $d_{1}$-dimensional
Wiener process on this space with complete, relative to
$\cF,P$, $\sigma$-fields $\cF_{t}$. Let $\sigma_{t},t\geq0$,
be a progressively measurable process with values in the set 
of $d\times d_{1}$-matrices and let $b_{t},t\geq0$, be an $\bR^{d}$-valued
progressively  measurable process.
Assume that for any $T\in[0,\infty)$ and $\omega$
\begin{equation}
                                             \label{8.19.2}
\int_{0}^{T}\big(|\sigma_{t}|^{2}+|b_{t}|)\,dt<\infty.
\end{equation}

Under this condition the stochastic process  
\begin{equation}
                                             \label{8.19.1}
x_{t}=\int_{0}^{t}\sigma_{s}\,dw_{s}
+\int_{0}^{t}b_{s}\,ds
\end{equation}
is well defined. Fix a nonnegative Borel $\gb$ on $\bR^{d}$.

\begin{assumption}
                                   \label{assumption 8.19.1} 
We have $\|\gb\| :=\|\gb\|_{L_{d}(\bR^{d})}<\infty$ and
\begin{equation}
                                             \label{8.19.3}
|b_{t}|\leq \gb(x_{t})\sqrt[d]{\det a_{t}}
\end{equation}
for  all $(\omega,t)$, where $a_{t}=(1/2)\sigma_{t}
\sigma^{*}_{t}$.

\end{assumption}

This assumption is supposed to hold
throughout the whole article.

Our goal is to establish
Aleksandrov's type estimates 
in the elliptic and parabolic settings and then apply them 
 to estimating the resolvent
operators of the corresponding elliptic
and parabolic operators in $L_{p}$
and $L_{p+1}$, respectively, where $p\geq d$.

Before stating our starting result proved in Section 
\ref{section 10.15.2} introduce
$$
B_{R}:=\{x\in\bR^{d}:|x|<R\}, 
$$
and let
$\tau_{R}(x)$ be the first exit time of $x+x_{t}$ from  
$B_{R}$ (equal to infinity if $x+x_{t}$ never exits from $B_{ R}$).
 Also let $\tau_{R}=\tau_{R}(0)$.
\begin{theorem}
                                     \label{theorem 10.13.1}
 
There is a constant $N_{d,\|\gb\|}$ (depending only on $d$ and 
$\|\gb\|$)
such that for any $R\in(0,\infty)$, $x\in\ B_{R} $,
and nonnegative Borel $f$ given on $\bR^{d}$  
we have
\begin{equation}
                                             \label{8.19.06}
E\int_{0}^{\tau_{ R}(x)}f(x+x_{t})\sqrt[d]{\det a_{t}}\,dt
\leq N_{d,\|\gb\|}R\|f\|_{L_{d}(B_{ R})}.
\end{equation}

\end{theorem}

As a result of future development of this theorem
we have the following proved in Section \ref{section 10.15.2}.
\begin{corollary}
                                       \label{corollary 10.20.2}
Suppose that for all $\omega\in\Omega$ and $t\geq0$ 
 we have $\tr a_{t}
\leq K $,
where $K\in(0,\infty)$ is a fixed number.
 Then
for $0\leq s <t<\infty$ and $R>0$ we have
\begin{equation}
                                             \label{12.4.2}
 P(\max_{r\in[s,t]}|x_{r}-x_{s}|\geq R\sqrt{t-s})
\leq 2e^{-\beta R^{2}/K},
\end{equation}
where $\beta=\beta(d,\|\gb\|)>0$.
In particular, for
any $n\geq0$ 
\begin{equation}
                                             \label{10.20.10} 
 E\max_{r\in[s,t]}|x_{r}-x_{s}|^{2n}\leq N (t-s)^{n},
\end{equation}
where $N=N(n,d,K,\|\gb\|)$.
\end{corollary}

\begin{remark}
                                           \label{remark 10.27.1}
The reader might have noticed that estimate
\eqref{10.20.10} has the same form for large $t-s$ and for small ones.
This is not the case if  $b$ is just bounded.
\end{remark}

The literature on the stochastic equations
with singular drift is quite impressive. We only mention one of the inspiring
ones \cite{Po_82} and some of
the most recent articles \cite{Hu_19},
\cite{Ro_16}, \cite{Zh_16}, \cite{Sh_13},
\cite{Fe_11},  \cite{Zh_05}. Also see  numerous
references therein. The goals in these articles
are much more ambitious than here, where we mostly
concentrate  around the results like Theorem
\ref{theorem 10.13.1} and Corollary
\ref{corollary 10.20.2}. In all of the above mentioned
sources apart from \cite{Ro_16} the coefficient
$b_{t}$ is allowed to depend on $x_{t}$ and $t$:
$b_{t}=b(t,x_{t})$ and it is supposed that
$$
\Big(\int\Big(\int|b(t,x)|^{p}\,dx\Big)^{q/p}\,dt\Big)^{1/q}
<\infty,
$$
where
\begin{equation}
                                           \label{11.13.1}
\frac{d}{p}+\frac{1}{q}<1.
\end{equation}
This excludes $p=d$ even
if $b(t,x)$ is independent of $t$.
In \cite{Ro_16} the setting is closer to ours
and again it is assumed that $p>d$.
Observe in passing that the parabolic PDE counterparts
of Theorem \ref{theorem 10.13.1} or of our
Theorem \ref{theorem 11.7.2} are obtained
in spaces with mixed norms in \cite{Na_15}
with $\leq$ in place of $<$ in \eqref{11.13.1}.
These results are not applicable in our situation,
in particular, because $x_{t}$ is not supposed to
solve any equation.

One may ask what happens if we only have $\gb\in L_{p}$
with $p<d$. It turns out that then estimate
\eqref{8.19.06} may break down. For instance,
let $dx_{t}=dw_{t}+b_{\varepsilon}(x_{t})\,dt$, $x_{0}=0$,
where $b_{\varepsilon}(x)=-(d/2)x|x|^{-2}I_{1>|x|>\varepsilon}$,
and in \eqref{8.19.06} take $x=0$, $f\equiv 1$, $R=1$.
Then simple computations in polar coordinates show
that the left-hand side of \eqref{8.19.06} tends to infinity
as $\varepsilon\downarrow 0$. At the same time, for any
$p<d$ the $L_{p}$-norms of $b_{\varepsilon}$ stay bounded.

On the basis of the results in this article,
the author intends to
show in a subsequent article 
that stochastic equations
with measurable diffusion and drift in $L_{d}$
admit solutions, that   $L_d$ in \eqref{8.19.06}
can be replaced with $L_p$ for some $p<d$
if the process $x_{t}$ is uniformly nondegenerate,
that It\^o's formula is applicable to
$u\in W^{2}_{p}$ for some $p<d$ and so on.
In particular, we will present some results concerning
estimates of the time spent in sets of small measure,
the probability to reach such sets, 
Fanghua Lin's estimates playing a major role in the Sobolev space theory
of {\em fully nonlinear elliptic equations\/},
boundary behavior of solutions of the corresponding
elliptic equations with first order coefficients in 
$L_{d}$, and the probability to pass through narrow tubes
(see \cite{Kr_19_1}).
It is worth saying that the main  driving force
behind our results is an   idea of Safonov
from \cite{Sa_10}.

Many results of this article could be
obtained by using the theory of fully nonlinear equations
presented,
for instance, in \cite{Kr_18}. This theory is 
ultimately based on Theorems \ref{theorem 8.18.1}
and \ref{theorem 11.7.4} and it turns out that
we only need these two theorem rather than rather 
sophisticated theory from \cite{Kr_18}
to derive our results.

The article is organized as follows.  
In Section \ref{section 10.15.2} we prove various 
extensions of Aleksandrov's type estimates for the elliptic
case (when $f$ in estimates like \eqref{8.19.06}
is independent of $t$) in bounded domains and in the whole space.
The latter allows us to also prove Corollary
\ref{corollary 10.20.2}.

In Section \ref{section 11.9.1} we give an application
of the results from Section \ref{section 10.15.2}
to showing that the resolvents of elliptic operators
with drift in $L_{d}$ decay in a normal way.
This is an improvement over some results in 
\cite{Kr_12}.

Section \ref{section 11.7.1} is devoted to
the parabolic case and estimates like
\eqref{8.19.06} but with $f(t,x+x_{t})$
in place of $f(x+x_{t})$. Its results are
applied in Section \ref{section 11.10.1}
to estimating the resolvent of parabolic
operators.

The final Section   \ref{section 10.11.1} is an
Appendix,
where we prove an auxiliary fact  related to
the theory of concave functions.

We finish the introduction with some notation and a stipulation
about constants. In the proofs of various results  we use
the symbol $N$  to denote finite 
nonnegative constants
which may change from one occurrence to another and,
 if in the statement of a result there are constants
called $N$ which are claimed to depend only on certain
parameters, then in the proof of the result
the constants $N$ also depend only on the same
parameters unless specifically stated otherwise.
Of course, if we write 
$
N=N(...),
$
 this means that $N$ depends only
on what is inside the parentheses. 

For $\delta\in(0,1)$ by $\bS_{\delta}$ we mean the set
of $d\times d$-symmetric matrices whose eigenvalues are
between $\delta$ and $\delta^{-1}$.
Introduce  
$$
B_{R}(x)=\{y:|y-x|<R\},\quad a_{\pm}=a^{\pm}=(1/2)(|a|\pm a),
$$
$$
D_{i}=\frac{\partial}{\partial x^{i}},
\quad D_{ij}=D_{i}D_{j} . 
$$
We use  the notation $u^{(\varepsilon)}=u*\zeta_{\varepsilon}$,
where
$\zeta_{\varepsilon}(x)=\varepsilon^{-d}\zeta(x/\varepsilon)$,
$\varepsilon>0$,
and $\zeta$ is a nonnegative $C^{\infty}$-function with
support in $B_{1}$ whose integral is equal to one.

\mysection{Elliptic case}

                                    \label{section 10.15.2}

Recall that Assumption \ref{assumption 8.19.1}
is supposed to be satisfied.

Our first result is Theorem \ref{theorem 10.13.1},
which is actually a particular case of Theorem 5.2
of \cite{Kr_86}, proved there by using a heavy artillery
from the theory of fully nonlinear elliptic 
partial differential equations. We provide its proof
based on the very initial knowledge of the solvability of the
Monge-Amp\`ere equations and estimates of its solutions,
which, actually, after a {\em long development\/} became also
one of the cornerstones of the theory 
of fully nonlinear elliptic 
partial differential equations.
The proof will be given after some preparations.
By setting $f=1$ in \eqref{8.19.06} we obtain the following.
\begin{corollary}
                                     \label{corollary 8.25.1}
We have
$$
E\int_{0}^{\tau_{R}(x)}\sqrt[d]{\det a_{t}}\,dt\leq NR^{2},
$$ 
where $N$ depends only on $d$ and $\|\gb\|$.

\end{corollary}

\begin{remark}
                             \label{remark 1.29.1}

The proof of Theorem \ref{theorem 10.13.1} consists of several steps.
First we let to the reader to check that the constant $N_{d,\|\gb\|}$
is independent of $R$ by replacing $x+x_{t}$ with $(x+x_{t})/R$.
In the new situation
condition \eqref{8.19.3} will change, but the new function $ \gb $
will have {\em the same\/} $L_{d}(\bR^{d})$-norm as the original one.
Therefore, it suffices to prove the theorem for $R=1$.
\end{remark}

Then the left-hand side of \eqref{8.19.06} will only increase
if we replace $\tau_{1}(x)$ with the first exit time
from $B_{2}(x)$. It follows that we may assume that
$x=0$ and $R=2$.

Now it becomes clear that
to prove Theorem \ref{theorem 10.13.1} it suffices
to prove the following.

\begin{theorem}
                                     \label{theorem 8.19.1}
 
There is a constant $N_{d,\|\gb\|}$
such that for any nonnegative Borel $f$ given on $B_{2 }$
we have
\begin{equation}
                                             \label{8.19.6}
E\int_{0}^{\tau }f(x_{t})\sqrt[d]{\det a_{t}}\,dt
\leq N_{d,\|\gb\|} \|f\|_{L_{d}(B_{2 })},
\end{equation}
where $\tau=\tau_{2}$.

\end{theorem}

To prove this theorem we need two lemmas in which
the assumptions of the theorem are supposed to hold
and we also need the following result proved in the Appendix.

\begin{theorem}
                                   \label{theorem 10.11.1}
Let $\theta$ be $0$ or $1$ and let $f $ be a nonnegative
function on $\bR^{d}$ such that $f^{d}$
 has finite integral over $B_{2}$ and $f=0$ outside $B_{2}$.
Then  there exists a convex function $z$ on $B_{4}$ which
is nonpositive  and such that

(a) we have in $B_{4}$
\begin{equation}
                                             \label{8.18.2}
|z|\leq 8\Phi\Big(\int_{B_{2}}f^{d}\,dx\Big),
\end{equation}
where $\Phi$ is the inverse function of
$$
\rho \to\int_{|p|\leq\rho}\frac{1}{(1+\theta|p|)^{d}}\, dp,
$$

(b) for any $\varepsilon\in(0,2)$
and nonnegative symmetric matrix $a$, in $B_{2 }$
we have
\begin{equation}
                                             \label{8.18.4}
a^{ij}D_{ij}z ^{(\varepsilon)} \geq
d\sqrt[d]{\det a}\big(f(1+\theta|Dz |)\big)^{(\varepsilon)}  .
\end{equation}

\end{theorem}
\begin{remark}
                                 \label{remark 8.19.1}
As is easy to see there is a constant $N_{d}$
depending only on $d$ such that
\begin{equation}
                                             \label{8.19.5}
8\Phi(t)\leq N_{d}t^{1/d}\quad\text{\rm if}
\quad \theta=0,\quad 8\Phi(t)\leq N_{d}t^{1/d}\exp(N_{d}t)
\quad\text{\rm if}
\quad \theta=1.
\end{equation}
\end{remark}

\begin{lemma}
                                           \label{lemma 8.19.1}
 
Let $\gamma$ be a stopping time such that $\gamma\leq\tau$
and  there is a constant $N$
such that \eqref{8.19.6}  
 holds with that $N$ and $\gamma$ in place of $\tau$
for any nonnegative Borel $f$ given on $B_{2}$.
Then it also holds with
($\gamma$ in place of $\tau$ and) an $N$ that
depends only on $d$ and $\|\gb\|$.

\end{lemma}

Proof. It suffices to concentrate on the $f$'s
that are continuous in $\bR^{d}$ and vanish
outside $B_{2}$.
Fix such an $f$ and define $z_{1}$ as a function
from Theorem \ref{theorem 10.11.1} when $\theta=0$.
By Remark \ref{remark 8.19.1}
$$
|z_{1}|\leq N_{d}F,\quad F:=\|f\|_{L_{d}(B_{2})}.
$$
Since, $z_{1}$ is convex and nonpositive, in $B_{3}$
we have
\begin{equation}
                                             \label{8.19.7}
|Dz_{1}(x)|\leq |z_{1}|/(4-|x|)\leq  |z_{1}|\leq  N_{d}F.
\end{equation}
Then define $z_{2}$ as a function
from Theorem  \ref{theorem 10.11.1} when $\theta=1$
and $f=\gb/d$. Finally, set $z=z_{1}+ N_{d}Fz_{2}$.

In light of \eqref{8.18.4} 
in $B_{2 }$, for any  nonnegative symmetric matrix $a$,
 we have
$$
a^{ij} D_{ij}z ^{(\varepsilon)} 
-\gb \sqrt[d]{\det a}|Dz^{(\varepsilon)} |\geq
d\sqrt[d]{\det a} f ^{(\varepsilon)} 
-\gb \sqrt[d]{\det a}|Dz_{1}^{(\varepsilon)} |
$$
$$
+ N_{d}F \sqrt[d]{\det a}\big(\gb(1+ |Dz_{2} |)\big)^{(\varepsilon)}
- N_{d}F\gb \sqrt[d]{\det a}|Dz_{2}^{(\varepsilon)} |.
$$
Owing to \eqref{8.19.7}, the last expression in $B_{2}$,
for $\varepsilon\in(0,1 )$, is greater than
$$
d\sqrt[d]{\det a} f ^{(\varepsilon)}- N_{d}F\sqrt[d]{\det a}
I_{\varepsilon},
$$
where
$$
I_{\varepsilon}:=\big(\gb(1+ |Dz_{2} |)\big)^{(\varepsilon)}
-\gb-\gb|Dz_{2}^{(\varepsilon)} |.
$$
Similarly to \eqref{8.19.7}, 
$I_{\varepsilon}\leq N(\gb+\gb^{(\varepsilon)})$, by 
the Lebesgue
theorem $I_{\varepsilon}\to0$ as $\varepsilon\downarrow0$
almost everywhere, and hence  as $\varepsilon\downarrow0$
we have
\begin{equation}
                                          \label{11.7.7}
 \|I_{\varepsilon} \|_{L_{d}(B_{2})}
\to 0 .
\end{equation}
By observing that thanks to \eqref{8.19.3}, for 
$t\leq\tau$, 
\begin{equation}
                                             \label{8.19.9}
a^{ij}_{t} D_{ij}z ^{(\varepsilon)} (x_{t})
+b^{i}_{t}Dz^{(\varepsilon)}(x_{t})  
\geq a^{ij} D_{ij}z ^{(\varepsilon)}(x_{t}) 
-\gb(x_{t}) \sqrt[d]{\det a_{t}}|Dz^{(\varepsilon)}(x_{t}) |
\end{equation}
and using It\^o's formula, we get that for small $\varepsilon$
and and $T\in[0,\infty)$
$$
Ez^{(\varepsilon)}(x_{\gamma \wedge T})
=z^{(\varepsilon)}(0)
+E\int_{0}^{\gamma \wedge T}\big(
a^{ij}_{t} D_{ij}z ^{(\varepsilon)} (x_{t})
+b^{i}_{t}Dz^{(\varepsilon)}(x_{t}) \big)\,dt
$$
$$
\geq z^{(\varepsilon)}(0)
+dE\int_{0}^{\gamma \wedge T}
\sqrt[d]{\det a_{t}}f^{(\varepsilon)}(x_{t})\,dt
$$
\begin{equation}
                                             \label{8.19.07}
-N_{d}FE\int_{0}^{\gamma}
\sqrt[d]{\det a_{t}}\,|I_{B_{2}}I_{\varepsilon}(x_{t})|\,dt.
\end{equation}

The last term in \eqref{8.19.07},
by absolute value, is less than a constant
independent of $\varepsilon$  times
the norm in \eqref{11.7.7} by assumption. Therefore,
it tends to zero as $\varepsilon\downarrow0$. The first term
in \eqref{8.19.07} is nonpositive. Therefore, by passing to the limit
in \eqref{8.19.07} first as $\varepsilon\downarrow0$ 
and then as $T\to\infty$ and also
using that $f$ is continuous we find
$$
E\int_{0}^{\gamma }f(x_{t})\sqrt[d]{\det a_{t}}\,dt
\leq d^{-1}|z(0)|
$$
$$
\leq d^{-1}N_{d}F+
d^{-1}N_{d}FN_{d}(\|\gb\|/d)^{1/d}\exp(N_{d}(\|\gb\|/d)^{d}).
$$
Hence \eqref{8.19.6} holds with $\gamma$ in place of $\tau$ and
$N=N_{d,\|\gb\|}$, where
\begin{equation}
                                             \label{8.19.8}
N_{d,\|\gb\|}:=d^{-1}N_{d} +
d^{-1}N_{d} N_{d}(\|\gb\|/d)^{1/d}\exp(N_{d}(\|\gb\|/d)^{d}).
\end{equation}
The lemma is proved.\qed

\begin{lemma}
                                           \label{lemma 8.19.2}
 
Let $\gamma$ be a stopping time such that $\gamma\leq\tau$
and    \eqref{8.19.6}  
 holds with   ($N_{d,\|\gb\|}$ from \eqref{8.19.8} and)
 $\gamma$ in place of $\tau$
for any nonnegative Borel $f$ given on $B_{2}$.
Introduce
$$
\gamma'=\tau\wedge\inf\{t\geq \gamma:
\int_{\gamma}^{t}|b_{s}|\,ds\geq 1\}.
$$
Then
there is a constant $N$
such that \eqref{8.19.6}  
 holds with that $N$ and $\gamma'$ in place of
$N_{d,\|\gb\|}$ and $\tau$, respectively,
for any nonnegative Borel $f$ given on $B_{2}$.

\end{lemma}

Proof. We take the same $f$, $z_{1}$, and $z_{2}$
as in the proof
of Lemma \ref{lemma 8.19.1} but use \eqref{8.19.9}
only for $t\leq \gamma$. On the interval
$[\gamma,\gamma')$ we use that  $|x_{t}|< 2 $, so that
(see \eqref{8.19.7})
 $|Dz^{(\varepsilon)}(x_{t})|\leq   dN_{d,\|\gb\|}F$
 implying  that
on $[\gamma,\gamma')$
$$
a^{ij}_{t} D_{ij}z ^{(\varepsilon)} (x_{t})
+b^{i}_{t}Dz^{(\varepsilon)}(x_{t}) 
\geq d\sqrt[d]{\det a_{t}}f^{(\varepsilon)}(x_{t})-
  dN_{d,\|\gb\|}F|b_{t}|.
$$
We also use that
$$
 E\int_{\gamma}^{\gamma'}|b_{t}|\,dt\leq 1.
$$
Hence, similarly to \eqref{8.19.07} we get
$$
Ez^{(\varepsilon)}(x_{\gamma'\wedge\tau_{\varepsilon}\wedge T})
=z^{(\varepsilon)}(0)
+E\int_{0}^{\gamma' \wedge T}\big(
a^{ij}_{t} D_{ij}z ^{(\varepsilon)} (x_{t})
+b^{i}_{t}Dz^{(\varepsilon)}(x_{t}) \big)\,dt
$$
$$
\geq z^{(\varepsilon)}(0)
+dE\int_{0}^{\gamma' \wedge T}
\sqrt[d]{\det a_{t}}f^{(\varepsilon)}(x_{t})\,dt
$$
\begin{equation}
                                             \label{8.19.10}
- N_{d}FE\int_{0}^{\gamma}
\sqrt[d]{\det a_{t}}|I_{\varepsilon}(x_{t})|\,dt
-dN_{d,\|\gb\|}F.
\end{equation}
As $\varepsilon\downarrow$ the term with $I_{\varepsilon}$
disappears in light of our assumption and of
what is said in the proof of Lemma  \ref{lemma 8.19.1}.
After that as in that proof we come to \eqref{8.19.6}
with   $N=2N_{d,\|\gb\|} $ and $\gamma'$ in place of 
$N_{d,\|\gb\|}$ and $\tau$, respectively.
The lemma is proved.\qed

{\bf Proof of Theorem \ref{theorem 8.19.1}}.
Set $\gamma_{0}=0$ and define recursively $\gamma_{n}$, $n=1,2,...$,
by
$$
\gamma_{n}=\tau\wedge\inf\{t\geq \gamma_{n-1}:
\int_{\gamma_{n-1}}^{t}|b_{s}|\,ds\geq  1 \}.
$$
In light of \eqref{8.19.2}, $\gamma_{n}\to\tau$ as $n\to\infty$.
Furthermore, since estimate \eqref{8.19.6} is obviously true
for $\gamma_{0}$ ($=0$) in place of $\tau$,   Lemma
\ref{lemma 8.19.2} implies that \eqref{8.19.6} 
is  true with a constant,
perhaps different from $N_{d,\|\gb\|}$ and with $\gamma_{1}$
in place of $\tau$. In that case Lemma \ref{lemma 8.19.1}
says that \eqref{8.19.6} is   true
for $\gamma_{1}$  in place of $\tau$. An obvious induction
proves that, for any $n$, \eqref{8.19.6} is   true
for $\gamma_{n}$  in place of $\tau$. Letting $n\to\infty$
we get the desired result. The theorem is proved.\qed

\begin{remark}
                                          \label{remark 8.21.1}
Theorem \ref{theorem 10.13.1} as many results below admits a natural
generalization with conditional expectations.
This generalization is obtained by tedious and not informative
repeating the proof with
obvious changes. We mean the following which we call
the {\em conditional\/} version of Theorem \ref{theorem 10.13.1}. Let $\gamma$ be a stopping time.
Then $x_{t}$ on the set $\{\gamma<\infty\}$ is given for $t\geq\gamma$
as
$$
x_{t}=x_{\gamma}+\int_{\gamma}^{t}\sigma_{t}\,dw_{t}
+\int_{\gamma}^{t}b_{t}\,dt.
$$
Also on the set $\{\gamma<\infty\}$
 introduce $\tau$ as the first time after $\gamma$
when $x_{t}$ exits from $B_{R}(x_{\gamma})$.  
Then   for any nonnegative Borel $f$ given on $\bR^{d}$
on the set $\{\gamma<\infty\}$ (a.s) we have
\begin{equation}
                                             \label{9.3.1}
E\Big[\int_{\gamma}^{\tau}f(x_{t})\sqrt[d]{\det a_{t}}\,dt
\mid \cF_{\gamma}\Big]
\leq N_{d,\|\gb\|}R \|f\|_{L_{d}(B_{R}(x_{\gamma}))}.
\end{equation}

Then standard measure-theoretic arguments show that
on the set $\{\gamma<\infty\}\times\bR^{d}$ there exists a function
$G_{\gamma,\tau}(x)=G_{\gamma,\tau}(\omega,x)$ which is
 $\cF_{\gamma}\otimes\cB(\bR^{d})$-measurable and
such that,
for any nonnegative Borel $f$ given on $\bR^{d}$,
on the set $\{\gamma<\infty\}$ (a.s) we have
$$
E\Big[\int_{\gamma}^{\tau}f(x_{t})\sqrt[d]{\det a_{t}}\,dt
\mid \cF_{\gamma}\Big]
=\int_{\bR^{d}}f(x)G_{\gamma,\tau}(x)\,dx.
$$
In addition \eqref{9.3.1} shows that on the set $\{\gamma<\infty\}$ (a.s.)
\begin{equation}
                                             \label{9.3.2}
\|G_{\gamma,\tau}\|_{L_{d/(d-1)}(B_{R}(x_{\gamma}))}
\leq N_{d,\|\gb\|}R   .
\end{equation}

Let us show briefly how to obtain \eqref{9.3.1}.
First by repeating the above proof one shows that (a.s.)
$$
E\Big[\int_{\gamma}^{\tau}f(x_{t}-x_{\gamma})\sqrt[d]{\det a_{t}}\,dt
\mid \cF_{\gamma}\Big]
\leq N_{d,\|\gb\|}R \|f\|_{L_{d}(B_{R} )}.
$$
Then by considering first continuous $f(x,y)$ and then
arbitrary nonnegative Borel, one gets that (a.s.)
$$
E\Big[\int_{\gamma}^{\tau}f(x_{t}-x_{\gamma},x_{\gamma})
\sqrt[d]{\det a_{t}}\,dt
\mid \cF_{\gamma}\Big]
\leq N_{d,\|\gb\|}R \|f(\cdot,x_{\gamma})\|_{L_{d}(B_{R} )}.
$$
Finally, the  substitution $f(x,y)=f(x+y)$ leads to
\eqref{9.3.1}.
\end{remark}

This remark allows us to generalize Corollary 
\ref{corollary 8.25.1}. For $0\leq s\leq t$ set
$$
\psi_{s,t}=
\int_{s}^{t}\sqrt[d]{\det a_{r}}\,dr,\quad
\psi_{t}=\psi_{0,t}.
$$

\begin{lemma}
                                     \label{lemma 10.20.1}
For $n=1,2,...$ and $x\in B_{R}$ we have
$$
I_{n}:=E\psi_{\tau_{R}(x)}^{n}
\leq N^{n}n!R^{2n},
$$ 
where $N$ depends only on $d$ and $\|\gb\|$.

\end{lemma}

Proof. Observe that
$$
I_{n+1}=(n+1)E\int_{0}^{\tau_{R}(x)}\sqrt[d]{\det a_{t}}\,
E\Big(\psi_{t,\tau_{R}(x)}^{n}
\mid \cF_{t}\Big)\,dt.
$$
If our assertion is true for a given $n$, then by its conditional
version
$$
I_{n+1}\leq (n+1)N^{n}n!R^{2n}
E\int_{0}^{\tau_{R}(x)}\sqrt[d]{\det a_{t}}\,dt.
$$
After that Corollary \ref{corollary 8.25.1} and the induction
on $n$ finish the proof.\qed

\begin{corollary} 
                                     \label{corollary 10.22.1}
There are  constants $N,\nu>0$ depending only on
$d$ and $\|\gb\|$ such that $E\exp(\nu
\psi_{\tau_{R}(x)}/R^{2})\leq N$ for any $R\in(0,\infty)$
and $x\in B_{R}$. In particular,
(by Chebyshev's inequality)
for any $t>0$,
$$
P(\psi_{\tau_{R}(x)}\geq t)\leq Ne^{-\nu t/R^{2}}.
$$
\end{corollary}

Corollary \ref{corollary 10.22.1} basically says that
$\tau_{R}(x)$ is smaller than a constant times $R^{2}$.
We want to show that in a sense the converse is also true.
For that, it is convenient to introduce
$$
\phi_{ t}=\phi_{0,t},\quad
\phi_{s,t}=\int_{s}^{t}\tr a_{r}\,dr.
$$

Here is a fact of immense importance for the future
development.

\begin{lemma}
                                     \label{lemma 8.21.1}
There is a constant $R=R_{d,\|\gb\|}\in[2,\infty)$ depending only on $d$ and 
$\|\gb\|$
such that 
\begin{equation}
                                             \label{8.21.2}
EI_{\tau<\infty}\exp\big(-\phi_{\tau}\big)
\leq 1/2,
\end{equation}
where $\tau=\tau_{R}$.

\end{lemma}

Proof. We use a Safonov's idea
from \cite{Sa_10}
based on the fact that
the integrals of $\gb^{d}$ should be small
an a variety of sets. Fix an integer $k\geq1$
 to be specified later
and define $\Gamma^{i}=B_{2(i+1)}\setminus \bar B_{2i}$.
Then for any
integer $n\geq k$ there are at least $m=n-k$ 
sets  $\Gamma^{i}$ with $i\leq n$ such that
$$
\|\gb\|_{L_{d}(\Gamma^{i})}\leq \|\gb\|/k^{1/d}.
$$
Let us call such sets ``good'' and let $i_{1}<...<i_{m}$
be such that each  $\Gamma^{i_{j}}$ is good. Introduce
$$
\tau^{j}=\inf\{t\geq0: |x_{t}|=2i_{j}+1\}\quad
(\inf\emptyset:=\infty),
$$
which (if finite) is the first time $x_{t}$ touches the middle sphere in 
$\Gamma^{i_{j}}$, and define $\gamma^{j}$ as the first
exit time of $x_{t}$ from  $B_{1}(x_{\tau^{j}})$ after 
$\tau^{j}$.

Observe that simple manipulations
show that
 the function $p (x):=\cosh  |x | $
satisfies
$$
a^{rs}D_{rs}p-p\,\tr a\leq0
$$
in $B_{1} $
for any symmetric nonnegative $d\times d$-matrix $a$.
By using It\^o's formula and observing
that $|Dp|\leq p\leq \cosh 1$ in $B_{1}$
 and  $|b^{r}_{t}D_{r}p(x_{t}-x_{\tau^{j}})|
\leq \gb(x_{t})\sqrt[d]{\det a_{t}}\cosh 1$ for
$t\in[\tau^{j},\gamma^{j})$ we see that for any $\varepsilon>0$  
on the set $\{\tau^{j}<\infty\}$ (a.s.)
$$
\cosh 1E\Big[ \exp\big(-\int_{\tau^{j}}^{\gamma^{j}}(\varepsilon+
\tr a_{t} )\,dt\big)
\mid \cF_{\tau^{j}}\Big]
$$
$$
=E\Big[ p(x_{\gamma^{j}}-
x_{\tau^{j}})\exp\big(-\int_{\tau^{j}}^{\gamma^{j}}(\varepsilon+
\tr a_{t} ) \,dt\big)
\mid \cF_{\tau^{j}}\Big]
$$
$$
\leq 1+\cosh 1E\Big[\int_{\tau^{j}}^{\gamma^{j}}
\gb(x_{t})\sqrt[d]{\det a_{t}}\,dt\mid \cF_{\tau^{j}}\Big]
$$
$$
= 1+\cosh 1E\Big[\int_{\tau^{j}}^{\gamma^{j}}
I_{\Gamma^{i_{j}}}\gb(x_{t})\sqrt[d]{\det a_{t}}
\,dt\mid \cF_{\tau^{j}}\Big].         
$$
By the way, the role of $\varepsilon$ is that the second term
in the above sequence makes sense even if $\gamma^{j}=\infty$.

By the choice of $\Gamma^{i_{j}}$ and 
Remark \ref{remark 8.21.1}
we can estimate the last term and conclude that
on the set $\{\tau^{j}<\infty\}$ (a.s.)
$$
\cosh 1E\Big[ \exp\big(-\int_{\tau^{j}}^{\gamma^{j}}(\varepsilon+
\tr a_{t} )\,dt\big)
\mid \cF_{\tau^{j}}\Big]\leq 1+(\cosh 1)N_{d,\|\gb\|}\|\gb\|/k^{1/d},
$$
\begin{equation}
                                             \label{8.21.1}
E\Big[ \exp\big(-\int_{\tau^{j}}^{\gamma^{j}}(\varepsilon+
\tr a_{t} )\,dt\big)\mid \cF_{\tau^{j}}\Big]\leq
 (\cosh 1)^{-1}+N_{d,\|\gb\|}\|\gb\|/k^{1/d}.
\end{equation}

We now choose $k$ and $n $, depending only on $d$ and $\|\gb\|$, so that
$$
\big((\cosh 1)^{-1}+N_{d,\|\gb\|}\|\gb\|/k^{1/d}\big)^{n-k}\leq 1/2
$$
and set
$$
R=2(n+1)\quad(\geq 2(i_{m}+1)).
$$

After that it only remains to observe that, obviously,
$$
E \exp\big(-\int_{0}^{\tau}(\varepsilon+\tr a_{t}) \,dt\big)
$$
is less than the expectation of the product
over $j=1,...,m$ of the left-hand sides
of \eqref{8.21.1}, and thus, less than $1/2$ for
any $\varepsilon>0$, and then send $\varepsilon\downarrow 0$.
The lemma is proved.\qed
\begin{corollary}
                                           \label{corollary 9.2.3}
For any   $R\in(0,\infty)$
\begin{equation}
                                             \label{9.2.4}
EI_{\tau_{R}<\infty}\exp\big(-(R_{d,\|\gb\|}/R)^{2}\phi_{\tau_{R}}\big)
\leq 1/2.
\end{equation}
  In particular, for $t_{0}:=R^{-2}_{d,\|\gb\|}\ln 3/2$
and any $R\in(0,\infty)$ we have
\begin{equation}
                                             \label{9.3.4}
P(\phi_{\tau_{R}}\leq t_{0}R^{2})\leq 3/4,
\end{equation}
provided that (a.s.)
\begin{equation}
                                             \label{9.3.40}
\int_{0}^{\infty}\tr a_{t}\,dt=\infty.
\end{equation}
\end{corollary}

Proof. Replacing $x_{t}$ with $(R_{d,\|\gb\|}/R)x_{t}$ shows that
it suffices to concentrate on $R=R_{d,\|\gb\|}$
(see Remark \ref{remark 1.29.1}). In that case 
\eqref{9.2.4} is identical with \eqref{8.21.2}. As far as
\eqref{9.3.4} is concerned it suffices to observe
that for $R=R_{d,\|\gb\|}$ the left-hand side equals
$$
P( \phi_{\tau_{R}}\leq \ln 3/2)=
P(\tau_{R}<\infty,\phi_{\tau_{R}}\leq \ln 3/2)
=P(\tau_{R}<\infty,e^{-\phi_{\tau_{R}}}\geq 2/3)
$$
$$
\leq (3/2)
EI_{\tau_{R}<\infty}e^{-\phi_{\tau_{R}}}\leq 3/4.
$$\qed

The following gives a more general form
to Corollary \ref{corollary 9.2.3}.
\begin{theorem}
                                      \label{theorem 10.16.01}
For any  $\kappa\in(0,1)$, $R\in(0,\infty)$, 
$x\in B_{\kappa R}$, and $\lambda\geq0$, 
\begin{equation}
                                             \label{11.3.1}
EI_{\tau_{R}(x)<\infty}\exp\big(-\lambda 
\phi_{\tau_{R}(x)}\big)
\leq 2e^{-\sqrt{\lambda} (1-\kappa)R/N},
\end{equation}
where $N=R_{d,\|\gb\|}/\ln2$. In particular, for  
  any $R,t\in(0,\infty)$ we have
\begin{equation}
                                             \label{10.16.2}
P(\phi_{\tau_{R}(x)}\leq  t R^{2})\leq 
2\exp\Big(-\frac{\beta(1-\kappa)^{2}}{t}\Big),
\end{equation}
where $\beta:=4^{-1}R^{-2}_{d,\|\gb\|}\ln^{2}2$ ($<1$), provided that 
\eqref{9.3.40} holds (a.s.).

\end{theorem}

Proof. Observe that $\phi_{\tau_{R}(x)}\geq \phi_{\tau'_{R}(x)}$,
where $\tau'_{R}(x)$ is the first exit time of $x+x_{t}$
from $B_{(1-\kappa)R}(x)$. It follows that in the proof of 
\eqref{11.3.1} we may assume that $\kappa=0$ and $x=0$.
Then, as usual (see Remark \ref{remark 1.29.1})
 we may assume that $R=1$. In that case
take $n\geq1$, to be specified later, and
introduce $\tau^{k}$, $k=1,...,n$, as the first exit
time of $x_{t}$ from   $B_{k/n}$. Also let
$\gamma^{k}$, $k=1,...,n$, be the first exit times
of $x_{t}$ from   $B_{1/n}(x_{\tau^{k-1}})$ after $\tau^{k-1}$
($\tau^{0}:=0$).
 Observe that, obviously $\phi_{\tau_{1}}
\geq \phi_{\tau^{0},\gamma^{1}}+...+\phi_{\tau^{n-1},\gamma^{n}}$
and by the conditional version of 
Corollary \ref{corollary 9.2.3} on the set $\{\tau^{k-1}<\infty\}$
(a.s.)
$$
E\Big(I_{\gamma^{k}<\infty}\exp\big(-R^{2}_{d,\|\gb\|}
n^{2}\phi_{\tau^{k-1},\gamma^{k}}\big)
\mid \cF_{\tau^{k-1}}\Big)\leq 1/2.
$$
After that   we get
$$
EI_{\tau_{1} <\infty}\exp\big(-R^{2}_{d,\|\gb\|} n ^{2}  
\phi_{\tau_{1} }\big)
\leq 1/2^{n},
$$
which for $n=\lfloor\sqrt{\lambda}/R_{d,\|\gb\|}\rfloor$
 yields \eqref{11.3.1}.

We prove \eqref{10.16.2} again assuming that $R=1$. We have
$$
P(\phi_{\tau_{1}(x)}\leq  t  )=P\big(
\tau_{1}(x)<\infty,\exp(-\lambda 
\phi_{\tau_{1}(x)}\geq \exp(-\lambda t)\big)
$$
$$
\leq 2 \exp( \lambda t-\sqrt{\lambda}(1-\kappa)/N).
$$
For $\sqrt{\lambda}=(1-\kappa)/(2Nt)$ we get \eqref{10.16.2}.
The theorem is proved. \qed

{\bf Proof of Corollary \ref{corollary 10.20.2}}.
First observe that \eqref{10.20.10} follows from \eqref{12.4.2}.
Next, by having in mind the
 conditional versions of our results   
we convince ourselves that while proving \eqref{12.4.2}
we may assume $s=0$. 
Finally, using scaling ($x_{t}\to cx_{t/c^{2}}$) reduces the 
case of general $t>0$  
to that of $t=1$. In that case, 
it only remains to
observe that
$$
P(\sup_{t\leq 1}|x_{t}|\geq R)=
P(\tau_{R} \leq 1)\leq P(\phi_{\tau_{R} }\leq K)
\leq 2e^{-\beta R^{2}/K}.
$$\qed

It is instructive to compare the following
with Corollary \ref{corollary 8.25.1}.
\begin{corollary}
                                       \label{corollary 7.29.1}
Suppose  
that \eqref{9.3.40} holds (a.s.). 
Then there is a constant $N=N(d, \|\gb\|)$
such that for any $R\in(0,\infty)$ 
\begin{equation}
                                                   \label{8.23.1}
NE\int_{0}^{\tau_{R}}\tr a_{t}\,dt\geq R^{2}.
\end{equation}

\end{corollary}

Indeed
$$
E\phi_{\tau_{R}}\geq  t_{0}R^{2}P(\phi_{\tau_{R}}>t_{0}R^{2})\geq 
(1/4)t_{0}R^{2}.
$$

Here is a technically convenient form of Corollary
\ref{corollary 7.29.1} to be used in a subsequent article
(see \cite{Kr_19_1}),
while estimating 
the time spend by $x_{t}$ in sets of small
measure.

\begin{corollary}  
                    \label{corollary 9.30.1}
Let $\delta\in(0,1)$ and assume that $\tr a_{t}\geq
\delta$ for all $(\omega,t)$. Then
for any $\kappa\in(0,1)$ there exists
$\nu=\nu(\kappa,d,\delta,\|\gb\|)>0$ such that for any $
R\in(0,\infty)$ and $x\in B_{\kappa R}$
\begin{equation}
                           \label{9.30.1}
E\tau\geq\nu R^{2},
\end{equation}
where $\tau=\tau_{R}(x)$.
\end{corollary}

Indeed,   $\tau\geq\gamma$,
where $\gamma$ is the first exit time of 
$x+x_{t}$ from $B_{(1-\kappa)R}(x)$ for which by Corollary
\ref{corollary 7.29.1}
$$
E\gamma\geq \delta^{-1}E\int_{0}^{\gamma}\tr a_{t} \,dt\geq
\delta^{-1}(1-\kappa)^{2}R^{2}\mu( d, \|\gb\|)
$$
with $\mu=\mu( d,\|\gb\|)>0$.

The following modestly looking version
of Corollary \ref{corollary 7.29.1} will
play a major role in a subsequent paper
while showing that $d$ in Theorems
\ref{theorem 10.13.1}  
can be  replaced with a smaller number
if $a_{t}\in\bS_{\delta}$ for all $\omega,t$
(see \cite{Kr_19_1}).

\begin{lemma}
                                             \label{lemma 11.7.2}
Let $a_{t}\in\bS_{\delta}$ for all $\omega,t$
and $\varepsilon\in(0,1]$.
Then there is a constant $N=N(d,\delta,\|\gb\|,\varepsilon)$ such that
for any $R\in(0,\infty)$ 
\begin{equation}
                                          \label{9.2.3}
NE\int_{0}^{\tau_{R}\wedge (\varepsilon R^{2})}e^{-t}\,dt
\geq R^{2}\wedge 1.
\end{equation}

\end{lemma}

Proof. Observe that, for $t\geq 0$, $e^{-t}\geq 1-t$, and
$$
I:=E\int_{0}^{\tau_{R}\wedge(\varepsilon R^{2})}e^{-t}\,dt
\geq E(\tau_{R}\wedge (\varepsilon R^{2}))-(1/2)
E(\tau_{R}\wedge (\varepsilon R^{2}))^{2}
$$
$$
\geq
E(\tau_{R}\wedge (\varepsilon R^{2}))-(1/2) R^{4}.
$$
Here by Theorem \ref{theorem 10.16.01}, for
$t\in (0,\varepsilon)$,
$$
E(\tau_{R}\wedge (\varepsilon R^{2}))\geq tR^{2}P(\tau_{R}>tR^{2})
\geq R^{2}t\Big(1-2\exp\Big(-\frac{\beta}{t}\Big)\Big),
$$
where $\beta=\beta(d,\delta,\|\gb\|)>0$,
so that for an appropriate $t=t(\beta,\varepsilon)>0$
$$
I\geq N^{-1}R^{2}-NR^{4},
$$
where the last expression is greater than $N^{-1}R^{2}$
for $R\leq R_{0}=R_{0}(d,\delta,\|\gb\|,\varepsilon)$ and \eqref{9.2.3}
holds for such $R$. However, for $R\geq R_{0}$
by Corollary \ref{corollary 9.2.3} for an appropriate
$s_{0}=s_{0}(d,\delta,\|\gb\|)>0$
$$
I\geq EI_{\tau_{R}>s_{0}R^{2}}
\int_{0}^{(s_{0}\wedge\varepsilon)R^{2}}e^{-t}\,dt\geq(1/4)
\int_{0}^{(s_{0}\wedge\varepsilon)R_{0}^{2}}e^{-t}\,dt.
$$
The lemma is proved. \qed

Theorem \ref{theorem 10.16.10}   
which is a simple
  corollary of Theorems \ref{theorem 10.13.1} and 
\ref{theorem 10.16.01}  can be used to prove
  Harnack's inequality for
diffusion processes with drift in $L_{d}$.

\begin{theorem}
                                     \label{theorem 10.16.10}
Assume that $a_{t}\in\bS_{\delta}$ for all $(\omega,t)$. Then
for any $\kappa\in(0,1)$ there is a 
function $q(\gamma)$, $\gamma\in(0,1)$,
depending only on $d,\delta,\|\gb\|,\kappa$, and, naturally, on $\gamma$,
such that for any $R\in(0,\infty)$, $x\in B_{\kappa R}$,
and closed $\Gamma\subset B_{R}$ satisfying
$|\Gamma|\geq \gamma|B_{R}|$ we have
$$
P(\tau_{\Gamma}(x)\leq \tau_{R}(x))\geq q(\gamma),
$$
where $\tau_{\Gamma}(x)$ is the first time the process $x+x_{t}$
hits $\Gamma$. Furthermore, $q(\gamma)\to 1$ as $\gamma\uparrow 1$.

\end{theorem}

Proof. By using scalings we reduce the general case
to the one in which $R=1$. In that case for any $\varepsilon>0$ we have
$$
P(\tau_{\Gamma}(x)> \tau_{1}(x))\leq
P\Big(\tau_{1}(x)=\int_{0}^{\tau_{1}(x)}I_{B_{1}\setminus\Gamma}(x+x_{t})
\,dt\Big)
$$
$$
\leq P(\tau_{1}(x)\leq\varepsilon  )+\varepsilon^{-1} 
E\int_{0}^{\tau_{1}(x)}I_{B_{1}\setminus\Gamma}(x+x_{t}) 
\,dt.
$$
In light of Theorems \ref{theorem 10.16.01} and 
\ref{theorem 10.13.1}
we can estimate the right-hand side and then obtain
$$
P(\tau_{\Gamma}(x)> \tau_{R}(x))\leq 2e^{-N/\varepsilon}
+N\varepsilon^{-1}|B_{1}\setminus\Gamma|^{1/d} 
$$
$$
\leq 2e^{-N/\varepsilon}   
+N\varepsilon^{-1}(1-\gamma)^{1/d}
$$
where the constants $N$  depend only on $d,\delta$, $\kappa$, and $\|\gb\|$.
By denoting
$$
q(\gamma)=1-\inf_{\varepsilon>0}\big(
2e^{-N/\varepsilon}
+N\varepsilon^{-1}(1-\gamma)^{1/d}\big),
$$
we get what we claimed. \qed

\begin{remark}
                                         \label{remark 10.20.1}
In a subsequent article (see \cite{Kr_19_1}) we will see
that   $Nq(\gamma)\geq \gamma^{\beta} $ for 
any $\gamma\in(0,1)$, where $N$ and $\beta>0$
depend only on $d,\delta,\|\gb\|,\kappa$.

\end{remark}

Next, we turn our attention to estimates
like in Theorem \ref{theorem 10.13.1}
 but on the infinite time interval.
\begin{lemma}
                                    \label{lemma 8.22.1}
 
For any Borel nonnegative $f$ vanishing outside
  $B_{1}$ and any $x\in\bR^{d}$ we have
\begin{equation}
                                       \label{11.9.3}
E\int_{0}^{\infty}e^{-\phi_{t}}f(x+x_{t})\sqrt[d]{\det a_{t}}\,dt 
  \le N 
\|f\| _{L_{d}(B_{1})},
\end{equation}
where $N$ depends only on $d$ and $\|\gb\|$. 

\end{lemma}

Proof. We may assume that $f$ is bounded. In  that case
introduce $\frM$ as the collection of all stopping times
and for $\tau\in\frM$ set
$$
\bar u_{\tau}:=\esssup_{\Omega} u_{\tau},\quad
\bar u=\sup_{\tau\in\frM}u_{\tau},
$$
where for a fixed $x\in\bR^{d}$
$$
u_{\tau}:=E\Big[
\int_{\tau}^{\infty}e^{-\phi_{\tau,t}}f(x+x_{t})
\sqrt[d]{\det a_{t}}\,dt
\mid \cF_{\tau}
\Big].
$$
Observe that, since $\sqrt[d]{\det a_{t}}\leq\tr a_{t}$ and 
$f$ is bounded, $\bar u<\infty$.

Take $R_{d,\|\gb\|}$ from Lemma \ref{lemma 8.21.1},
 take a $\tau\in\frM$, and define $\gamma$ as the first
time after $\tau$ when $|x_{t}-x_{\tau}|\geq R_{d,\|\gb\|}$.
On the set where $\tau<\infty$ and $|x+x_{\tau}|\geq 1+R_{d,\|\gb\|}$
we have (recall that $f$ is zero outside $B_{1}$)
$$
u_{\tau}=E\Big[e^{-\phi_{\tau,\gamma}}E
\Big[\int_{\gamma}^{\infty}e^{-\phi_{\gamma,t}}
f(x+x_{t})\sqrt[d]{\det a_{t}}\,dt\mid\cF_{\gamma}\Big]
\mid \cF_{\tau}\Big]
$$
$$
=E \Big[e^{-\phi_{\tau,\gamma}}u_{\gamma}\mid \cF_{\tau}\Big]\leq 
\bar u
E \Big[e^{-\phi_{\tau,\gamma}} \mid \cF_{\tau}\Big],
$$
where by the conditional version of Lemma \ref{lemma 8.21.1}
the last conditional expectation is less than $1/2$
(a.s.).

On the set where $\tau<\infty$ and $|x+x_{\tau}|\leq 1+R_{d,\|\gb\|}$
we have $B_{1}\subset B_{\rho}(x+x_{\tau})$, where
$\rho=3+2R_{d,\|\gb\|}$. The choice of $\rho$ is dictated
by the fact that the outside of $B_{\rho}(x+x_{\tau})$
is also outside of $B_{1+R_{d,\|\gb\|}}$. Therefore,
if we define $\theta$ as the first time after $\tau$
when $|x_{t}-x_{\tau}|\geq\rho$, then by the above
$u_{\theta}\leq \bar u/2$ (a.s.). After that
we observe that
$$
u_{\tau}=E\Big[
\int_{\tau}^{\theta}e^{-\phi_{\tau,t}}f(x+x_{t})\sqrt[d]{\det a_{t}}\,dt
+e^{-\phi_{\tau,\theta}}u_{\theta}\mid \cF_{\tau}
\Big]
$$
$$
\leq E\Big[
\int_{\tau}^{\theta}e^{-\phi_{\tau,t}}f(x+x_{t})\sqrt[d]{\det a_{t}}\,dt
\mid \cF_{\tau}
\Big]+\bar u/2,
$$
where by the conditional version of Theorem 
\ref{theorem 10.13.1} the last conditional expectation
is less than a constant $N$ depending only on $d$ and $\|\gb\|$
times the $L_{d}$-norm of $f$.

Thus in both cases of the location of $x+x_{\tau}$ we have (a.s.)
$$
u_{\tau}\leq N\|f\|_{L_{d}(B_{1})}+(1/2)\bar u.
$$
The arbitrariness of $\tau$ leads to the conclusion
that 
$$
\bar u\leq N\|f\|_{L_{d}(B_{1})}+(1/2)\bar u
$$
and the result follows.\qed

\begin{corollary}
                               \label{corollary 8.22.2}
For any Borel nonnegative $f$ vanishing outside
  $B_{1}$ and $x\in\bR^{d}$ we have
$$
I:=E\int_{0}^{\infty}e^{-\phi_{t}}f(x+x_{t})\sqrt[d]{\det a_{t}}\,dt 
  \le N e^{-\mu|x|}
\|f\| _{L_{d}(B_{1})},
$$
where $N$ and $\mu>0$ depend  only on $d$ and $\|\gb\|$.

\end{corollary}

Indeed, we only need to consider the case that
$|x|\geq 2$ in which we let $\tau$ be the first 
exit time of $x_{t}$ from $B_{|x|-1}$
and observe that in the notation from Lemma
\ref{lemma 8.22.1}
$$
I=EI_{\tau<\infty}e^{-\phi_{\tau}} u_{\tau}
\leq N\|f\|_{L_{d}(B_{1})}
EI_{\tau<\infty}e^{-\phi_{\tau }}.
$$ 
 After that it only remains to observe
that $EI_{\tau<\infty}e^{-\phi_{\tau}}\leq Ne^{-|x|/N}$
by Theorem \ref{theorem 10.16.01}.

\begin{theorem}
                                        \label{theorem 8.22.1}
Let   $p\geq d$.
 Then there exists   constants $N$ and $\mu>0$,
depending only on $d,p$, and $\|\gb\|$, such that for any
$\lambda>0$ and Borel nonnegative $f$ given on $\bR^{d}$
we have
\begin{equation}
                                                   \label{7.29.2}
u:=E\int_{0}^{\infty}e^{-\lambda \phi_{t}}
f(x_{t})\sqrt[d]{\det a_{t}}\,dt   \le N\lambda^{d/(2p)-1}
\|\Psi_{\lambda}^{-1} f\| _{L_{p}(\bR^{d})},
\end{equation}
where $\Psi_{\lambda}(x)=\exp( \sqrt{\lambda}\mu|x|)$.

\end{theorem}

Proof. Since the case of arbitrary $\lambda>0$
is reduced to the case $\lambda=1$ by simply replacing $x_{t}$
with $x_{t}\sqrt{\lambda}$, we only concentrate on
$\lambda=1$.
Take a nonnegative 
$\zeta\in C^{\infty}_{0}(B_{1})$ which integrates to one,
 for $y\in\bR^{d}$ set $f_{y}(x)=\zeta(x-y)f(x)$,
and introduce
$$
u(y)=E\int_{0}^{\infty}e^{-  \phi_{t}}
f_{y}(x_{t})\sqrt[d]{\det a_{t}}\,dt .
$$
By Corollary \ref{corollary 8.22.2}, for any  $y\in\bR^{d}$,
$$
u(y)\leq Ne^{-2\mu|y|}\|f_{y}\|_{L_{d}(\bR^{d})}
$$
$$
\leq Ne^{- \mu|y|}\|\Psi_{1}^{-1} f_{y}\|_{L_{d}(\bR^{d})}
\leq Ne^{-\mu|y| }\|\Psi_{1}^{-1} f_{y}\|_{L_{p}(\bR^{d})}.
$$
After that it only remains to note that
$$
\int_{\bR^{d}}u(y)\,dy=u,\quad
\int_{\bR^{d}}\|\Psi_{1}^{-1} f_{y}\|^{p}_{L_{p}(\bR^{d})}\,dy
=\|\zeta\|_{L_{p}(\bR^{d})}\|\Psi_{1}^{-1} f\|_{L_{p}(\bR^{d})}
$$
and use H\"older's inequality ($q=p/(p-1)$)
$$
\int_{\bR^{d}}e^{-\mu|y|}\|\Psi_{1}^{-1} f_{y}\|_{L_{p}(\bR^{d})}\,dy
\leq\Big[\int_{\bR^{d}}e^{-q\mu|y|}\,dy\Big]^{1/q}
\Big[\int_{\bR^{d} }\|\Psi_{1}^{-1} f_{y}\|^{p}_{L_{p}(\bR^{d})}\,dy\Big]^{1/p}.
$$
The theorem is proved. \qed

Here is somewhat unexpected consequence of
Theorem \ref{theorem 8.22.1}.

{\bf The second proof of
\eqref{10.20.10} in Corollary \ref{corollary 10.20.2}}. 
As in the proof of Corollary \ref{corollary 10.20.2}
given after Theorem \ref{theorem 10.16.01}  
 we may assume $s=0$ and $t=1$. Using H\"older's
inequality allows us to concentrate on $n\geq1$. 
In this case 
observe that  using It\^o's formula
easily yields that
$$
|x_{1}|^{2n}\leq N\int_{0}^{1}(|x_{s}|^{2n-1}\gb(x_{s})+|x_{s}|^{2(n-1)})
\,ds+2n\int_{0}^{1}|x_{s}|^{2(n-1)}x^{*}_{s}\sigma_{s}\,dw_{s},
$$
where in light of \eqref{7.29.2} with  $p= d $  
 the last term is
a square integrable martingale, so that the expectation of its 
supremum over $r\leq 1$ is dominated by 
a constant times the square root of
$$
I:=E\int_{0}^{1}|x_{s}|^{4n-2}\,ds\leq e 
E\int_{0}^{1}e^{- s}|x_{s}|^{4n-2}\,ds.
$$
 We use \eqref{7.29.2} with  
an appropriate $\lambda=\lambda(d,\delta)$   to get that
$$
I\leq N \|\Psi^{-1}_{\lambda}|\cdot|^{4n-2}
\|_{L_{d}(\bR^{d})}=N .
$$
 Similarly,
$$
E\int_{0}^{1} |x_{s}|^{2(n-1)} 
\,ds\leq N \|\Psi^{-1}_{\lambda}|\cdot|^{2(n-1)}
\|_{L_{p}(\bR^{d})}=N .
$$

Finally,
$$
E\int_{0}^{1} |x_{s}|^{2n-1}\gb(x_{s})\,ds
\leq N \|\Psi^{-1}_{\lambda}|\cdot|^{2n-1}\gb
\|_{L_{d}(\bR^{d})}\leq N,
$$
where the last inequality holds because $\Psi^{-1}_{\lambda}
|x|^{2n-1}$ is bounded. 
This proves \eqref{10.20.10}. \qed

We finish the section with one more result which will be extended
in a subsequent article to wider range of $p$
but more narrow set of processes (see \cite{Kr_19_1}).

\begin{theorem}
                                        \label{theorem 11.23.1}
Let   $p\geq d$.
 Then there exists   constants $N$ and $\mu>0$,
depending only on $d,p$, and $\|\gb\|$,
and $R_{0}=R_{0}(d,\|\gb\|)$, such that for any
$\lambda>0$, $R\in[0,\infty)$,
 and Borel nonnegative $f$ given on $\bR^{d}$
we have
$$
E\int_{0}^{\infty}e^{-\lambda \phi _{t}(B_{R}^{c})}
f(x_{t})\sqrt[d]{\det a_{t}}\,dt  
$$
\begin{equation}
                                                   \label{11.23.1}
 \le N (R\sqrt{\lambda}+R_{0}) ^{2-d/p} 
\lambda^{d/(2p)-1} 
\|\Psi_{R,\lambda}^{-1} f\| _{L_{p}(\bR^{d})},
\end{equation}
where $\Psi _{R,\lambda}(x)=\exp\big(\sqrt{\lambda}\mu
\dist(x,B_{R+R_{0}/\sqrt{\lambda}})\big)$ and
$$
\phi _{t}(B_{R}^{c})=\int_{0}^{t}\tr a_{s}I_{|x_{s}|\geq R}\,ds.
$$
\end{theorem}

Proof. The change of coordinates $x_{t}=\sqrt{\lambda}y_{t}$
shows that we may assume that $\lambda=1$.
For $R=0$ the result follows
from Theorem \ref{theorem 8.22.1}.
In the general case, take $R_{0}=R(d,\|\gb\|)=R_{d,\|\gb\|}  $ 
from Lemma \ref{lemma 8.21.1} and define recursively, $\gamma^{0}=0$  
$$
\tau^{n}=\inf\{t\geq   \gamma^{n}:x_{t}\not\in B_{R+R_{0}}\},
\quad \gamma^{n+1}
=\inf\{t\geq   \tau^{n}:x_{t}\in \bar B_{R }\}.
$$
Also for simplicity write $\phi^{R}_{t}$
in place of $\phi _{t}(B_{R}^{c})$.

By the conditional version  of Theorem \ref{theorem 8.22.1}
   on the set ${\tau^{n}<\infty}$ (a.s.)
$$
 E\Big(
\int_{\tau^{n}}^{\gamma^{n+1}}e^{-  \phi^{R}_{t}}
f(x_{t})\sqrt[d]{\det a_{t}}\,dt\mid \cF_{\tau^{n}}\Big)
$$
$$
=e^{- \phi^{R}_{\tau^{n}}} E\Big(
\int_{\tau^{n}}^{\gamma^{n+1}}e^{-  (\phi_{t}-
\phi_{\tau^{n}})}
f(x_{t})\sqrt[d]{\det a_{t}}\,dt\mid \cF_{\tau^{n}}\Big)
$$
$$
\leq e^{- \phi_{\tau^{R}_{n}}}
 N 
\|\Psi_{R,1}^{-1} f\| _{L_{p}(\bR^{d})}.
$$
Furthermore, by the choice of $R_{0}$ and the conditional
version of Lemma \ref{lemma 8.21.1} for $n\geq 1$
$$
EI_{\tau^{n}<\infty}e^{-  \phi^{R}_{\tau^{n}}}
\leq EI_{\tau^{n-1}<\infty}
 e^{- \phi^{R}_{\tau^{n-1}}}E\Big(I_{\gamma^{n}<\infty}
e^{-  (\phi_{\gamma^{n}}-\phi_{\tau^{n-1}})}
\mid \cF_{\tau^{n-1}}\Big)
$$
$$
\leq (1/2)EI_{\tau^{n-1}<\infty}
 e^{-  \phi^{R}_{\tau^{n-1}}}.
$$
It follows that for $n\geq0$
\begin{equation}
                                                   \label{11.24.1}
 EI_{\tau^{n}<\infty}e^{- \phi^{R}_{\tau^{n}}}
\leq (1/2)^{n},
\end{equation}
\begin{equation}
                                                   \label{11.23.2}
 E 
\int_{\tau^{n}}^{\gamma^{n+1}}e^{- \phi^{R}_{t}}
f(x_{t})\sqrt[d]{\det a_{t}}\,dt 
\leq N(1/2)^{n} 
\|\Psi_{R,1}^{-1} f\| _{L_{p}(\bR^{d})}.
\end{equation}

On the other hand, by Theorem \ref{theorem 10.13.1}
and H\"older's inequality
on the set $\{\gamma^{n}<\infty\}$ (a.s.)
$$
e^{ \phi^{R}_{\gamma^{n}}} E\Big(
\int_{\gamma^{n}}^{\tau^{n }}e^{-  \phi^{R}_{t}}
f(x_{t})\sqrt[d]{\det a_{t}}\,dt\mid \cF_{\gamma^{n}}\Big) 
$$
$$
\leq E\Big(
\int_{\gamma^{n}}^{\tau^{n }} 
f(x_{t})\sqrt[d]{\det a_{t}}\,dt\mid \cF_{\gamma^{n}}\Big)
$$
$$
\leq N(R+R_{0}) 
  \| f\| _{L_{d}(B_{R+R_{0} )})}\leq N
(R+R_{0}) ^{2-d/p}
  \| f\| _{L_{p}(B_{R+R_{0} )})}
$$
$$
=N
(R+R_{0}) ^{2-d/p}
  \| \Psi_{R,\lambda}^{-1} f\| _{L_{p}(B_{R+R_{0} )})},
$$
where the last equality holds because $\Psi_{R,\lambda}=1$
on $B_{R+R_{0}}$. Furthermore, for $n\geq 1$,  
$$
EI_{\gamma^{n}<\infty}e^{-  \phi^{R}_{\gamma^{n}}}
\leq EI_{\tau^{n-1}<\infty}e^{-  \phi^{R}_{\tau^{n-1}}},
$$
so that owing to \eqref{11.24.1}
$$
EI_{\gamma^{n}<\infty}e^{-  \phi^{R}_{\gamma^{n}}}
\leq(1/2)^{n-1},
$$
which is also true for $n=0$.
Hence for $n\geq0$
$$
 E 
\int_{\gamma^{n}}^{\tau^{n }}e^{-  \phi^{R}_{t}}
f(x_{t})\sqrt[d]{\det a_{t}}\,dt 
\leq N(1/2)^{n-1}(R+R_{0}) ^{2-d/p}
  \|\Psi_{R,1}^{-1} f\| _{L_{p}(B_{R+R_{0} )})}.
$$
After that it only remains to take into account
\eqref{11.23.2} and the fact that ($R_{0}\geq2$ and)
$$
\int_{0}^{\infty}=\sum_{n=0}^{\infty}\Big(
\int_{\gamma^{n}}^{\tau^{n}}+\int_{\tau^{n}}^{\gamma^{n+1}}\Big)
$$
The theorem is proved.         \qed
 
 \mysection{An application to elliptic equations}
 
                                    \label{section 11.9.1}
Let $a(x)$ be a Borel measurable function
on $\bR^{d}$ with values in $\bS_{\delta}$, where $\delta
\in(0,1)$ is a fixed constant. Let $b(x)$
 be a Borel measurable function
on $\bR^{d}$ with values in $\bR^{d}$ 
such that 
$$
\|b\|_{L_{d}(\bR^{d})}\leq \|b\|,
$$
where $\|b\|<\infty$ is a fixed number. Define
$$
L=(1/2)a^{ij}D_{ij}+b^{i}D_{i}\quad \big(D_{i}=\frac{\partial}{\partial x^{i}},
\quad D_{ij}=D_{i}D_{j}\big).
$$ 

Here is the result of this section.
Such results play a crucial role in the proof
that one can pass to the limit under the sign 
of  fully nonlinear elliptic operators
when the arguments of these operators converge only in
a very weak sense (see, for instance,
Section 4.2 in \cite{Kr_18}).

Theorem \ref{theorem 10.14.2} for
 $R=\infty$ is obtained in \cite{Kr_12}
however
with $\lambda$ in \eqref{10.14.4} 
restricted from below by a constant
depending on how fast $\|(|b|-\mu)_{+}\|
_{L_{d}(\bR^{d})}\to 0$ as $\mu\to\infty$.
This is caused, in part, by the fact that $b=b_{1}+b_{2}$,
where $b_{1}$ is bounded and $b_{2}\in L_{d }$,
in \cite{Kr_12}.
\begin{theorem}
                                   \label{theorem 10.14.2}
Let $p\geq d $ and $R\in(0,\infty]$. Then there exists
a constant
$N=N(d,\delta,\|b\|)\geq0$ such that for any $\lambda>0$
and $u\in W^{2}_{p,\loc}(B_{R})\cap C(\bar B_{R})$ 
($B_{\infty}=\bR^{d}$, $C(\bR^{d})$
is the set of bounded continuous functions on $\bR^{d}$)
 we have
\begin{equation}
                                           \label{10.14.4}
\lambda\|u_{+}\|_{L_{p}(B_{R/2})}
\leq N\|(\lambda u-Lu)_{+}\|_{L_{p}(B_{R})}
+N\lambda R^{d/p}e^{-R\sqrt{\lambda}/N}\sup_{\partial B_{R}}u_{+},
\end{equation}
where the last term should be dropped if $R=\infty$.
\end{theorem}

Proof.  First we note that scaling
the coordinates allows us to assume that $\lambda=1$. 
By having in mind the possibility to approximate
$B_{R}$ from inside by similar domains,
we see that we may assume that 
$u\in W^{ 2}_{p }(B_{R})$ and $R<\infty$.
We may also assume that the norm on the right in
\eqref{10.14.4} is finite. Then we can replace $L$
with $L_{n}:=I_{|b|\geq n}\Delta+I_{|b|< n}L$
and if \eqref{10.14.4} is true for $L_{n}$
in place of $L$ we can pass to the limit
by the dominated convergence theorem.
Thus, we may assume that $b$ is bounded.
In that case it suffices to prove
\eqref{10.14.4} for $u\in C^{ 2}(  \bar B_{R})$.
Finally, by using mollifications we reduce
this case further to the case of smooth
bounded $a,b$.

In that case take $x_{0}\in B_{R/2}$, and
denote by $x_{t}$ a solution of 
$$
x_{t}=x_{0}+\int_{0}^{t}\sqrt{ a(x_{s})}\,dw_{s}
+\int_{0}^{t}b(x_{s}) \,ds
$$
existing on a probability space carrying a $d$-dimensional
Wiener process $w_{t}$.
 By It\^o's formula we conclude that
for any $T\in(0,\infty)$  
$$
u(x_{0})=E\int_{0}^{T\wedge\tau_{R}}e^{-t}f(x_{t})\,dt
+Ee^{-T\wedge\tau_{ R }} u(x_{T\wedge\tau_{ R }}),
$$
where $\tau_{R}$ is the first exit time of $x_{t}$
from $B_{R}$.
By sending $T\to\infty$ we obtain
\begin{equation}
                                           \label{10.15.3}
u(x_{0})=E\int_{0}^{\tau_{ R }}e^{-t}f(x_{t})\,dt
+Ee^{- \tau_{ R }} u(x_{ \tau_{ R }})=:I (x_{0})
+J (x_{0}).
\end{equation}

Note that, owing to \eqref{11.3.1},   
 $
Ee^{- \tau_{ R}}\leq 2e^{-R/N}
 $,
so that
$$
J_{+} (x_{0})\leq 2e^{-R/N}\sup_{\partial B_{R}}u_{+},
\quad 
 \|J_{+} \|_{L_{p}(B_{R/2})}
\leq  
 N R^{d/p}e^{-R/N}\sup_{\partial B_{R}}u_{+}.
$$

To estimate $I (x_{0})$,
we use the same method as in the proof of  
Theorem \ref{theorem 8.22.1}, first setting $f=0$
outside $B_{R}$ and observing that
$$
I (x_{0})\leq E\int_{0}^{\infty}e^{-t}f(x_{t})\,dt
=:\bar I (x_{0}).
$$
Take a nonnegative 
$\zeta\in C^{\infty}_{0}(B_{1})$ which integrates to one,
 for $y\in\bR^{d}$ set $f_{y}(x)=\zeta(x-y)f(x)$,
 introduce
$$
I(x_{0},y)=E\int_{0}^{\infty}e^{- t}
(f_{y}(x_{t}))_{+} \,dt 
$$
and observe that by Corollary \ref{corollary 8.22.2}
$$
 I(x_{0},y) \leq Ne^{-\nu|x_{0}-y|}\|(f_{y})_{+}\|_{L_{p}(\bR^{d})}.
$$
Also obviously
$$
\bar I(x_{0})\leq\int_{\bR^{d}}I(x_{0},y)\,dy,
$$
so that
$$
\bar  I_{+}(x_{0}) \leq N\int_{\bR^{d}}
e^{-\nu|x_{0}-y|}\|(f_{y})_{+}\|_{L_{p}(\bR^{d})}\,dy,
$$
where the integral is the convolution of two
functions. Hence
$$
\|\bar I_{+}\|_{L_{p}(\bR^{d})}\leq
N\int_{ \bR^{d} }e^{-\nu| y|}\,dy\Big( \int_{ \bR^{d} }
\Big(\int_{ \bR^{d} }\zeta^{p}(x-y)|f_{+}(x)|^{p}\,dx\Big)\,dy
\Big)^{1/p}
$$
$$
=N\|f_{+}\|_{L_{p}(\bR^{d})}
=N\|f_{+}\|_{L_{p}(R_{R})}
$$
and the theorem is proved.  \qed

\mysection{Parabolic  case}

                                    \label{section 11.7.1}

Recall that Assumption \ref{assumption 8.19.1}
is supposed to be satisfied.
Define
$$
C_{T,R}(t_{0},x_{0})=[t_{0},t_{0}+T)\times B_{R}(x_{0}),\quad 
C_{ R}(x_{0})=[0,\infty)\times B_{R}(x_{0}),
$$
$$
C_{T,R} =C_{T,R}(0,0),\quad
C_{ R} =C_{ R}(0),\quad \partial_{t}=\frac{\partial}{\partial t}.
$$

Here is our first result.
\begin{theorem}
                                     \label{theorem 11.7.2}
 
There is a constant $N_{d,\|\gb\|}$ (depending only on $d$ and 
$\|\gb\|$)
such that for any $R\in(0,\infty)$, $x\in\ B_{R}$,
and nonnegative Borel $f$ given on $C_{R}$  
we have
\begin{equation}
                                             \label{11.7.3}
E\int_{0}^{\tau_{ R}(x)}f(t,x+x_{t})\sqrt[d+1]{\det a_{t}}\,dt
\leq N_{d,\|\gb\|}R^{d/(d+1)}\|f\|_{L_{d+1}(C_{ R})}.
\end{equation}

\end{theorem}

As in the case of Theorem \ref{theorem 10.13.1}
to prove Theorem \ref{theorem 11.7.2} it suffices
to prove the following.

\begin{theorem}
                                     \label{theorem 11.7.3}
 
There is a constant $N_{d,\|\gb\|}$
such that for any nonnegative Borel $f$ given on $C_{2 }$
we have
\begin{equation}
                                             \label{11.7.4}
E\int_{0}^{\tau }f(t,x_{t})\sqrt[d+1]{\det a_{t}}\,dt
\leq N_{d,\|\gb\|} \|f\|_{L_{d+1}(C_{2 })},
\end{equation}
where $\tau =\tau_{2}$.

\end{theorem}

We need the following analytic result 
that is a particular case of Theorem 2 in \cite{Kr_76}
in which we use the notation
$f^{(\varepsilon)}=f*\xi_{\varepsilon}$, where
$\xi_{\varepsilon}(t,x)=\varepsilon^{-d-1}
\xi(t/\varepsilon,x/\varepsilon)$ and $\xi$
is a nonnegative $C^{\infty}$ function with unit integral
and support in $(-1,0)\times B_{1}$.

\begin{theorem}
                                     \label{theorem 11.7.4}
 
Let a nonnegative function
$f\in L_{d+1}(C_{4}) $ be such that $f=0$ outside $C_{2}$.
Then on $C_{4}$ there exists a bounded nonpositive
function $z(t,x)$ such that

a) $z(t,x)$ is convex in $x$ for $x\in B_{4}$,
$t\geq0$, and increasing in $t$ for $x\in B_{4}$;

b) for all nonnegative symmetric $d\times d$ matrices 
$a=(a^{ij})$ and $\varepsilon\in(0,2)$ we have on 
$C_{2}$ that
\begin{equation}
                                               \label{11.7.1}
\partial_{t}z^{(\varepsilon)}
+a^{ij}D_{ij}z^{(\varepsilon)} \geq
\alpha \sqrt[d+1]{\det a} f^{(\varepsilon)} ,
\end{equation}
where $\alpha=\alpha(d)>0$ is a constant;

c) there exists a constant $N_{d}=N(d)$ such that
\begin{equation}
                                               \label{11.7.2}
|z|\leq N_{d}\|f\|_{L_{d+1}(C_{ 2})}.
\end{equation}
\end{theorem}

{\bf Proof of Theorem \ref{theorem 11.7.3}}.
It suffices to concentrate on continuous $f$'s as in 
Theorem \ref{theorem 11.7.4}.
Take such an $f$ and take $ z $ from
 Theorem \ref{theorem 11.7.4}. Set $z_{1}= z$, denote
$$
F=\|f\|_{L_{d+1}(C_{ 2})},
$$
and with $z_{2}$ as in the proof of Lemma \ref{lemma 8.19.1}
and $N_{d}$ from \eqref{11.7.2}
introduce $z=z_{1}+ N_{d}Fz_{2}$.  

In light of \eqref{11.7.1} for all small $\varepsilon>0$
and nonnegative symmetric 
$d\times d$ matrix $a$
in
$C_{2}$   
 we have
$$
\partial_{t}z ^{(\varepsilon)}+a^{ij} D_{ij}z ^{(\varepsilon)} 
-\gb \sqrt[d]{\det a}|Dz^{(\varepsilon)} |\geq
\alpha\sqrt[d+1]{\det a} f ^{(\varepsilon)} 
-\gb \sqrt[d]{\det a}|Dz_{1}^{(\varepsilon)} |
$$
$$
+ N_{d}F \sqrt[d]{\det a}\big(\gb(1+ |Dz_{2} |)\big)^{(\varepsilon)}
- N_{d}F\gb \sqrt[d]{\det a}|Dz_{2}^{(\varepsilon)} |.
$$
Owing to \eqref{8.19.7},  the last expression 
on $C_{2}$ is greater than
$$
\alpha\sqrt[d+1]{\det a}
 f ^{(\varepsilon)}- N_{d}F\sqrt[d]{\det a}
I_{\varepsilon},
$$
where
$$
I_{\varepsilon}:=\big(\gb(1+ |Dz_{2} |)\big)^{(\varepsilon)}
-\gb-\gb|Dz_{2}^{(\varepsilon)} |
$$
is a function of $x$ alone, for which
\eqref{11.7.7} holds. 

By observing that thanks to \eqref{8.19.3}, for 
$t\leq\tau $.
$$
b^{i}_{t}Dz^{(\varepsilon)}(t,x_{t})  
\geq   
-\gb(x_{t}) \sqrt[d ]{\det a_{t}}\,|Dz^{(\varepsilon)}(t,x_{t}) |
$$
and using It\^o's formula, we get that for small $\varepsilon$
and $T\in(0,\infty)$
$$
0\geq Ez^{(\varepsilon)}( 
\tau\wedge T ,x_{ 
\tau \wedge T})
$$
$$
=z^{(\varepsilon)}(0)
+E\int_{0}^{ \tau \wedge T}\big(
\partial_{t}z ^{(\varepsilon)} (t,x_{t})+
a^{ij}_{t} D_{ij}z ^{(\varepsilon)} (t,x_{t})
+b^{i}_{t}Dz^{(\varepsilon)}(t,x_{t}) \big)\,dt
$$
$$
\geq z^{(\varepsilon)}(0)
+\alpha E\int_{0}^{ \tau \wedge T}  
\sqrt[d+1]{\det a_{t}}\,f^{(\varepsilon)}(t,x_{t})\,dt
$$
\begin{equation}
                                             \label{11.7.8}
- N_{d}FE\int_{0}^{\tau }
\sqrt[d]{\det a_{t}}\,|I_{\varepsilon}(x_{t})|\,dt.
\end{equation}

The last term in \eqref{11.7.8}  
 tends to zero as $\varepsilon\downarrow0$
in light of \eqref{11.7.7} and Theorem
\ref{theorem 10.13.1}.   By passing to the limit
in \eqref{11.7.8}   as $\varepsilon\downarrow0$ 
and $T\to\infty$
and also
using that $f$ is continuous we find
$$
E\int_{0}^{\tau}f(t,x_{t})\sqrt[d+1]{\det a_{t}}\,dt
\leq \alpha^{-1}|z(0)|
$$
$$
\leq \alpha^{-1}N_{d}F+
\alpha^{-1}N_{d}FN_{d}(\|\gb\|/d)^{1/d}\exp(N_{d}(\|\gb\|/d)^{d}).
$$
The theorem is proved. \qed

Below, for $\Gamma\subset \bR^{d+1}$, by $|\Gamma|$
we mean its Lebesgue measure.
The next theorem will
be used in a subsequent article
to prove the Harnack inequalities
for harmonic and caloric functions associated
with diffusion processes with drift in $L_{d}$
(see \cite{Kr_20}).
\begin{theorem}
                                     \label{theorem 11.8.1}
Assume that $a_{t}\in\bS_{\delta}$ for all $(\omega,t)$. Then
for any $\kappa\in(0,1)$ there is a 
function $q(\gamma)$, $\gamma\in(0,1)$,
depending only on $d,\delta,\|\gb\|,\kappa$, and, naturally, on $\gamma$,
such that for any $R\in(0,\infty)$, $x\in B_{\kappa R}$,
and closed $\Gamma\subset C_{R^{2},R}$ satisfying
$|\Gamma|\geq \gamma|C_{R^{2},R}|$ we have
$$
P(\tau_{\Gamma}(x)\leq \tau_{R^{2},R}(x))\geq q(\gamma),
$$
where $\tau_{\Gamma}(x)$ is the first time the process 
$(t,x+x_{t})$
hits $\Gamma$ and $\tau_{R^{2},R}(x)$ is its first exit time from
$C_{R^{2},R}$. Furthermore, 
$q(\gamma)\to 1$ as $\gamma\uparrow 1$.

\end{theorem}

Proof. By using scalings we reduce the general case
to the one in which $R=1$. In that case for any $\varepsilon\in(0,1)$ we have
$$
P(\tau_{\Gamma}(x)> \tau_{1,1}(x))\leq
P\Big(\tau_{1,1}(x)=\int_{0}^{\tau_{1,1}(x)}I_{C_{1,1}
\setminus\Gamma}(t,x+x_{t})
\,dt\Big)
$$
$$
\leq P(\tau_{1 }(x)\leq\varepsilon  )+\varepsilon^{-1} 
E\int_{0}^{\tau_{1,1}(x)}I_{C_{1,1}\setminus\Gamma}(t,x+x_{t}) 
\,dt.
$$
In light of Theorems \ref{theorem 10.16.01} and 
\ref{theorem 11.7.2}
we can estimate the right-hand side and then obtain
$$
P(\tau_{\Gamma}(x)> \tau_{1,1}(x))\leq 2e^{-N/\varepsilon}
+N\varepsilon^{-1}|C_{1,1}\setminus\Gamma|^{1/(d+1)} 
$$
$$
\leq 2e^{-N/\varepsilon}   
+N\varepsilon^{-1}(1-\gamma)^{1/(d+1)}
$$
where the constants $N$  depend only on $d,\delta$, $\kappa$, and $\|\gb\|$.
By denoting
$$
q(\gamma)=1-\inf_{\varepsilon\in(0,1)}\big(
2e^{-N/\varepsilon}
+N\varepsilon^{-1}(1-\gamma)^{1/(d+1)}\big),
$$
we get what we claimed. \qed

Next, we turn our attention to estimates
like in Theorem \ref{theorem 11.7.2}
 but on the infinite time interval.
\begin{lemma}
                                    \label{lemma 11.9.2}
 
For any $ x \in  \bR^{d}$, 
and Borel nonnegative $f$ on $\bR^{d+1}$
vanishing outside
  $C_{1 } $   we have
\begin{equation}
                                       \label{11.9.4}
E\int_{0}^{\infty}e^{-\phi_{t}}
f(t,x+x_{t})\sqrt[d+1]{\det a_{t}}\,dt 
  \le N 
\|f\| _{L_{d+1}(C_{1 } )},
\end{equation}
where $N$ depends only on $d$ and $\|\gb\|$. 

\end{lemma}

Proof. Following the proof of
Lemma \ref{lemma 8.22.1} we may assume that $f$ is bounded. In  that case
introduce $\frM$ as the collection of all stopping times
and for $\tau\in\frM$ set
$$
\bar u_{\tau}:=\esssup_{\Omega} u_{\tau},\quad
\bar u=\sup_{\tau\in\frM}u_{\tau},
$$
where,  for a fixed $\varepsilon>0$,
$$
u_{\tau}:=E\Big[
\int_{\tau}^{\infty}e^{-\varepsilon(t-\tau)-\phi_{\tau,t}}
f(t, x+x_{t})
\sqrt[d+1]{\det a_{t}}\,dt
\mid \cF_{\tau}
\Big].
$$

Observe that since $\sqrt[d+1]{\varepsilon\det a_{t}}
\leq\varepsilon+\tr a_{t}$ and 
$f$ is bounded, $\bar u<\infty$. This allows us to repeat
literally the proof of Lemma \ref{lemma 8.22.1}
with only one change that instead of
Theorem \ref{theorem 10.13.1} one should refer
to Theorem \ref{theorem 11.7.2}.
Then we get \eqref{11.9.4} with $\varepsilon t+\phi_{t}$
in place of $\phi_{t}$. After that it only remains to set
$\varepsilon\downarrow 0$. \qed

In the same way as Corollary \ref{corollary 8.22.2}
is obtained we arrive at the following.

\begin{corollary}
                                  \label{corollary 11.10.1}
For any Borel nonnegative $f$ vanishing outside
  $C_{1}$ and $x\in\bR^{d}$ we have
$$
E\int_{0}^{\infty}e^{-\phi_{t}}f(t,x+x_{t})
\sqrt[d+1]{\det a_{t}}\,dt 
  \le N e^{-\mu|x|}
\|f\| _{L_{d+1}(C_{1})},
$$
where $N$ and $\mu>0$ depend  only on $d$ and $\|\gb\|$.

\end{corollary}

Next, by repeating almost literally
the proof of
Theorem \ref{theorem 8.22.1} in case $p=d$
we obtain the following.

\begin{theorem}
                                        \label{theorem 11.10.1}
 
 There exists   constants $N$ and $\mu>0$,
depending only on $d $  and $\|\gb\|$, such that for any
$\lambda>0$ and Borel nonnegative $f$ given on $\bR^{d+1}$
we have
\begin{equation}
                                                   \label{7.29.02}
 E\int_{0}^{\infty}e^{-\lambda \phi_{t}}
f(t,x_{t})\sqrt[d+1]{\det a_{t}}\,dt   \le N
\lambda^{-d/(2d+2)}
\|\Psi_{\lambda}^{-1} f\| _{L_{d+1}(\bR^{d+1})},
\end{equation}
where $\Psi_{\lambda}(x)=\exp( \sqrt{\lambda}\mu|x|)$.

\end{theorem} 

By plugging $e^{-\lambda t}f$ in place of $f$
in \eqref{7.29.02} and then using H\"older's
inequality on its right-hand side we get
the following.
\begin{corollary}
                                   \label{corolary 11.10.3}
Let $p\geq d+1$.
There exists   constants $N$ and $\mu>0$,
depending only on $d,p$, and $\|\gb\|$, such that for any
$\lambda>0$ and Borel nonnegative $f$ given on $\bR^{d+1}$
we have
\begin{equation}
                                               \label{11.10.5}
 E\int_{0}^{\infty}e^{-\lambda t-\lambda \phi_{t}}
f(t,x_{t})\sqrt[d+1]{\det a_{t}}\,dt   \le N
\lambda^{ (d+2)/(2p)-1}
\|\Phi_{\lambda}  f\| _{L_{p}(\bR^{d+1})},
\end{equation}
where $\Phi_{\lambda}(x)=\exp( -\sqrt{\lambda}\mu|x|
-\lambda t/2)$.

\end{corollary}

\mysection{An application to parabolic equations}
 
                                    \label{section 11.10.1}
Let $a(t,x)$ be a Borel measurable function
on $\bR^{d+1}$ with values in $\bS_{\delta}$, where $\delta
\in(0,1)$ is a fixed constant. Let $b(t,x)$
 be a Borel measurable function
on $\bR^{d+1}$ with values in $\bR^{d}$ 
such that there exists a function $\gb\in L_{d}(\bR^{d})$  
for which
$$
|b(t,x)|\leq \gb(x) 
$$
for all $t$ and $x$. Set
$$
\|\gb\|=\|\gb\|_{L_{d}(\bR^{d})},
$$
and define
$$
L=\partial_{t}+(1/2)a^{ij}D_{ij}+b^{i}D_{i}.
$$ 
Also for $t\in[0,R^{2})$ set
$$
 C_{R^{2},R}(t)=[t,R^{2})\times B_{R},\quad
\partial'C_{R^{2},R}(t)=\partial C_{R^{2},R}(t)\setminus(\{t\}
\times B_{R}),
$$
and for $t=0$ drop the argument $t$ in the above
notation.

Here is the result of this section.
Such results play a crucial role in the proof
that one can pass to the limit under the sign 
of  fully nonlinear elliptic operators
when the arguments of these operators converge only in
a very weak sense (see, for instance,
Section 4.2 in \cite{Kr_18}).

Theorem \ref{theorem 11.10.3} for
 $R=\infty$ is obtained in \cite{Kr_12}
however
with $\lambda$ in \eqref{11.10.4} 
restricted from below by a constant
depending on how fast $\|(|b|-\mu)_{+}\|
_{L_{d+1}(\bR^{d+1})}\to 0$ as $\mu\to\infty$.
This is caused, in part, by the fact that $b=b_{1}+b_{2}$,
where $b_{1}$ is bounded and $b_{2}\in L_{d+1}$,
in \cite{Kr_12}.
\begin{theorem}
                                   \label{theorem 11.10.3}
Let $p\geq d+1 $ and $R\in(0,\infty]$. Then there exists
a constant
$N=N(d,\delta,\|\gb\|)\geq0$ such that for any $\lambda>0$
and $u\in W^{1,2}_{p,\loc}(C_{R^{2},R})\cap C(\bar C_{R^{2},R})$ 
($C_{\infty,\infty}=\{(t,x):t\geq0\}$, $C(C_{\infty,\infty})$
is the set of bounded continuous functions on 
$C_{\infty,\infty}$)
 we have
$$
\|u(t_{0},\cdot)\|_{L_{p}(B_{R/2})}
\leq N\lambda^{-(p-1)/p}\|(\lambda u-Lu)_{+}\|_{L_{p}(C_{R^{2},R}(t_{0}))}
$$
\begin{equation}
                                         \label{11.11.1}
+NR^{d/p}e^{-R\sqrt{\lambda}/N}\sup_{\partial'C_{R^{2},R}(t_{0})}u_{+} 
\end{equation}
for any $t_{0}\in[0,R^{2}/4]$ and 
$$
\lambda\|u_{+}\|_{L_{p}(C_{R^{2}/4,R/2})}
\leq N\|(\lambda u-Lu)_{+}\|_{L_{p}(C_{R^{2},R})}
$$
\begin{equation}
                                           \label{11.10.4} 
+N\lambda R^{(d+2)/p}e^{-R\sqrt{\lambda}/N}\sup_{\partial'
 C_{R^{2},R}}u_{+},
\end{equation}
where in both estimates 
the last terms should be dropped if $R=\infty$.
\end{theorem}

Proof. We follow the proof of Theorem 
\ref{theorem 10.14.2} and convince ourselves
that it suffices to concentrate on smooth
$a,b$, and $u$. Also it suffices to prove
\eqref{11.11.1} for $\lambda=1$ and
\eqref{11.10.4} for $\lambda=2$.

We start with
\eqref{11.11.1}, take $(t_{0},x_{0})\in C_{R^{2}/4,R/2}$, and
denote by $x_{t}$, $t\geq t_{0}$, a solution of 
$$
x_{t}=x_{0}+\int_{t_{0}}^{t}\sqrt{ a(s,x_{s})}\,dw_{s}
+\int_{t_{0}}^{t}b(s,x_{s}) \,ds
$$
existing on a probability space carrying a $d$-dimensional
Wiener process $w_{t}$. Denote by $\tau$
the first time $t$ when $(t,x_{t})$ exits from $C_{R^{2} ,R }$
after time $t_{0}$ that is the minimum of $R^{2}-t_{0}$
($\geq (3/4)R^{2}$)
and the first time $\gamma$ when
 $x_{t}$ exits from $B_{R}$ after time $t_{0}$.
Define $f=u-Lu$, $\bar f=fI_{C_{R^{2},R} } $.

 By It\^o's formula we conclude that
\begin{equation}
                                           \label{11.10.7}
u(t_{0},x_{0})=E\int_{t_{0}}^{\tau}
e^{-(t-t_{0})}\bar f(t, x_{t})\,dt
+E 
u(\tau ,x_{\tau })
e^{-(\tau-t_{0})}=:I+J.
\end{equation}

Note that, for any $\varepsilon\in (0,1)$
 $$
Ee^{-(\tau-t_{0})}\leq P(\tau-t_{0}< \varepsilon(R^{2}-t_{0}))
+e^{-(3/4)\varepsilon R^{2}},
 $$
where the first term on the right is $P(\gamma-t_{0}<
\varepsilon(R^{2}-t_{0}))$ which owing to \eqref{10.16.2}
is dominated by twice $\exp(-\beta/\varepsilon)$,
where $\beta=\beta(d,\delta,\|\gb\|)>0$.
It follows that, if $R>1$, by taking $\varepsilon=1/R$,
we get
$$
Ee^{-(\tau-t_{0})}\leq 2e^{-\beta R}
$$
with another $\beta=\beta(d,\delta,\|\gb\|)>0$,
for which we may safely assume that $\beta\leq 1$.
However, if $R<1$, 
$$
Ee^{-(\tau-t_{0})}\leq e^{\beta}e^{-\beta R}.
$$
It follows that
$$
J\leq e e^{-\beta R}\sup_{\partial'C_{R^{2},R}(t_{0})}u_{+}.
$$
 
To estimate $I $,
we use the same method as in the proof of  
Theorem \ref{theorem 8.22.1}.
Take a nonnegative 
$\zeta\in C^{\infty}_{0}(B_{1})$ which integrates to one,
 for $y\in\bR^{d}$ set $\bar f_{y}(t,x)=\zeta(x-y)\bar f(t,x)$,
 introduce
$$
I( y)=E\int_{t_{0}}^{\infty}e^{- (t-t_{0})}
(\bar f_{y}(t,x_{t}))_{+} \,dt 
$$
and observe that by Corollary \ref{corollary 11.10.1}
$$
 I(y) \leq Ne^{-\mu|x_{0}-y|}\|(\bar f_{y})_{+}\|_{L_{p}(\bR^{d+1})}.
$$

As a result,
\begin{equation}
                                              \label{11.10.9}
u_{+}(t_{0},x_{0})\leq N
\int_{\bR^{d}}e^{-\mu|x_{0}-y|}
\|(\bar f_{y})_{+}\|_{L_{p}(\bR^{d+1})}\,dy+
e e^{-\beta R}\sup_{\partial'C_{R^{2},R}(t_{0})}u_{+}.
\end{equation}
Here the first term on the right is a convolution
whose $L_{p}(\bR^{d})$-norm is less than the $L_{1}$-norm
of $e^{-\mu|x|}$ times 
$$
\Big(\int_{\bR^{d}} 
\|(\bar f_{y})_{+}\|^{p}_{L_{p}(\bR^{d+1})}\,dy\Big)^{1/p}
=N\|\bar f_{+}\|_{L_{p}(\bR^{d+1})}=
N\|f_{+}\|_{L_{p}(C_{R^{2},R}(t_{0}))}.
$$
After that \eqref{11.11.1} for $\lambda=1$ follows.

To prove \eqref{11.10.4} for $\lambda=2$, observe that
by taking $e^{-(t-t_{0})}u(t,x)$
in place of $u$ in \eqref{11.11.1} with $\lambda=1$
we obtain
$$
\|u(t_{0},\cdot)\|^{p}_{L_{p}(B_{R/2})}
\leq N\int_{0}^{R^{2}}I_{t_{0}\leq t}e^{-p(t-t_{0})}
\|(2u-Lu)_{+}(t,\cdot)\|^{p}_{L_{p}(B_{R })}\,dt
$$
$$
+N R^{d}e^{-pR /N}\sup_{\partial'
 C_{R^{2},R}}(u_{+})^{p}.
$$
By integrating through this estimate with respect 
to $t_{0}\in[0,R^{2}/4]$ we immediately obtain
\eqref{11.10.4} with $\lambda=2$.
The theorem is proved.  \qed

\begin{remark}
                                     \label{remark 11.11.1}
In light of Theorem \ref{theorem 11.7.2},
estimates \eqref{11.11.1}  and \eqref{11.11.1}
are of little use when $R$ is bounded and $\lambda$
is small. The real strength of this estimates reveals
when $\lambda$ is large or $R=\infty$.
\end{remark}

\mysection{Appendix: Proof of Theorem  
\protect\ref{theorem 10.11.1} }
                                          \label{section 10.11.1}

Let $z(x)$ be a finite convex function on $B_{4}$.
For each point $x_{0}\in B_{4}$ define $\phi(z,x_{0})$ 
as the collection
of $p\in\bR^{d}$ for each of which there exists $b\in\bR$
such that the plane $z=(p,x)+b$ is a supporting plane
for the graph of $z(x)$ at point $(z(x_{0}),x_{0})$.
For any Borel set $\Gamma\subset B_{4}$ Aleksandrov 
\cite{Al_58} sets 
$$
\phi(z,\Gamma)
=\bigcup_{x\in \Gamma}\phi(z,x),\quad \nu(z,\Gamma)
=\Vol \phi(z,\Gamma)
$$
 and proves that $\nu(z,\Gamma)$ is a measure
finite on any compact subset of $B_{4}$. He calls
$\phi(z,\Gamma)$ the normal image of $\Gamma$ 
relative to $z$.

We are going to consider the equation
\begin{equation}
                                             \label{8.18.1}
\int_{\phi(z,\Gamma)}\frac{1}{(1+\theta|p|)^{d}}\, dp 
=\int_{\Gamma}f^{d}(x)\,dx,
\end{equation}
which is an equation about the unknown
 convex $z$ which should be satisfied
for any Borel $\Gamma\subset B_{4}$, where $\theta$ is a parameter
and $f$ is a given function.

The following is a particular case of Theorem 4
of \cite{Al_58}, with \eqref{8.18.2} being a particular case
of estimate (2,8) of \cite{Al_58}.

\begin{theorem}
                                   \label{theorem 8.18.1}
Let $\theta$ be $0$ or $1$ and let $f $ be a nonnegative
function on $B_{4}$ such that $f^{d}$ has finite integral over $B_{4}$.
Then  equation \eqref{8.18.1} has a solution $z$ which
is nonpositive in $B_{4}$ and satisfies  estimate \eqref{8.18.2}.

\end{theorem}

According to \cite{Al_58} the left-hand side 
of \eqref{8.18.1} is a measure for any convex $z$.
Recall (see, for instance, \cite{Kr_82}) that for any convex $z$ on $B_{4}$
its gradient $Dz$ is well defined almost everywhere
and is a monotone function. It admits an extension
to the whole of $B_{4}$ and the extension denoted
by $(z_{1},...,z_{d})$ is still 
monotone. According to Theorem 2.1 of \cite{Kr_82}
the first generalized derivatives of the extension
are (signed) measures, finite on any compact subset
of $B_{4}$. Denote by $z^{(0) }_{i,j}(x)$ the density
of the measure $D_{j}z_{i}(dx)$

By Theorem 3.2 of \cite{Kr_82}
the density of the absolutely continuous part
of the left-hand side of \eqref{8.18.1} for any
convex $z$ equals
$$
\frac{1}{(1+\theta|Dz(x)|)^{d}}\det \big(z^{(0) }_{i,j}(x)\big)
$$

Next, recall that
$\zeta\in C^{\infty}_{0}(B_{1})$
is  a nonnegative radially symmetric function  which integrates to one,
  for $\varepsilon>0$ we define $\zeta_{\varepsilon}(x)
=\varepsilon^{-d}\zeta(x/\varepsilon)$, and for any distribution
$u$ we use the notation $u^{(\varepsilon)}=u*\zeta_{\varepsilon}$.
One knows that if $u$ is a measure, then $u^{(\varepsilon)}
\to u^{(0)}$ as $\varepsilon\downarrow 0$ almost everywhere,
where $u^{(0)}$ is the density of the absolutely continuous part
of the measure $u$. Owing to this and the fact that the
$D_{i}z$'s are the  generalized derivatives of $z$
coinciding with $z_{i}$'s (as generalized functions),
we get (a.e.)
$$
z^{(0) }_{i,j}=\lim_{\varepsilon\downarrow0}
\big(D_{j}z_{i}\big)^{(\varepsilon)}
=\lim_{\varepsilon\downarrow0}
D_{j}\big(z_{i}\big)^{(\varepsilon)}
=\lim_{\varepsilon\downarrow0}
D_{j}\big(D_{i}z\big)^{(\varepsilon)}
$$
$$
=\lim_{\varepsilon\downarrow0}
\big(D_{ij} z\big) ^{(\varepsilon)}=z_{ij}^{(0)},
$$
where $z_{ij}^{(0)}$ is
the density of the absolutely continuous part
of the measure $D_{ij}z$ which is the generalized derivative
of order $D_{ij}$ of $z$. Just in case, recall that
these  generalized derivatives are signed measures
existing for any convex $z$.
In this way we arrive at the first assertion of the 
following corollary, which is stated without proof in a slightly more
general situation by Aleksandrov
 in \cite{Al_58} and claimed to be a simple consequence
of his arguments in \cite{Al_39}.

\begin{corollary}
                                       \label{corollary 8.18.1}
The solution $z$ from Theorem \ref{theorem 8.18.1}
satisfies (a.e. $B_{4}$)
\begin{equation}
                                             \label{8.18.3}
\det(z_{ij}^{(0)})=f^{d}(1+\theta|Dz(x)|^{d}).
\end{equation}
Furthermore, for any $\varepsilon\in(0,2)$
and nonnegative symmetric matrix $a$,
\eqref{8.18.4} holds in $B_{2}$.  

\end{corollary}

To prove the second assertion of the corollary, observe that,
since the matrix of generalized derivatives $(D_{ij}z)$
is nonnegative,   owing to the inequality between
the arithmetic and the geometric means,
we have that in the sense of generalized
functions
$$
\sqrt[d]{\det(z_{ij}^{(0)})}\sqrt[d]{\det a}\,dx
\leq (1/d)a^{ij} z_{ij}^{(0)}\,dx
\leq(1/d) a^{ij}D_{ij}z(dx).
$$
Hence,
$$
(1/d) a_{ij}D_{ij}z(dx)\geq
\sqrt[d]{\det a} (f(1+\theta|Dz |)\,dx
$$
and it only remains to take convolutions of both sides
with $\zeta_{\varepsilon}$.\qed

\end{document}